\documentclass[11pt,leqno]{siamltex}
\usepackage[top=1in, bottom=1in, left=1in, right=1in]{geometry}
\usepackage{multirow}
\usepackage{xcolor}

\usepackage{amsmath,amsfonts,amssymb,amscd,mathrsfs}
\usepackage{extarrows,textcomp}
\usepackage{graphicx,subfigure,epstopdf}
\usepackage{diagbox,multirow}
\usepackage{algorithm,algorithmic}
\usepackage{caption}
\usepackage{bm}
\usepackage{float}
\usepackage{lineno}
\usepackage{xcolor}
\usepackage{makeidx,hyperref}

\newcommand{\td}{\text{d}}

\title{Structure Probing Neural Network Deflation}

\author{Yiqi Gu
\vspace{0.1in}\\
Department of Mathematics, National University of Singapore, 10 Lower Kent Ridge Road, Singapore, 119076 ({\tt matguy@nus.edu.sg})
\vspace{0.1in}\\
Chunmei Wang
\vspace{0.1in}\\
Department of Mathematics \& Statistics, Texas Tech University, 1108 Memorial Circle, Lubbock, TX 79409, USA ({\tt chunmei.wang@ttu.edu})
  \vspace{0.1in}\\
Haizhao Yang
  \vspace{0.1in}\\
  Department of Mathematics, Purdue University, 150 N University St, West Lafayette, IN 47907,  USA ({\tt haizhao@purdue.edu})
}

\begin{document}
\maketitle
\begin{abstract}
Deep learning is a powerful tool for solving nonlinear differential equations, but usually, only the solution corresponding to the flattest local minimizer can be found due to the implicit regularization of stochastic gradient descent. This paper proposes a network-based structure probing deflation method to make deep learning capable of identifying multiple solutions that are ubiquitous and important in nonlinear physical models. First, we introduce deflation operators built with known solutions to make known solutions no longer local minimizers of the optimization energy landscape. Second, to facilitate the convergence to the desired local minimizer, a structure probing technique is proposed to obtain an initial guess close to the desired local minimizer. Together with neural network structures carefully designed in this paper, the new regularized optimization can converge to new solutions efficiently. Due to the mesh-free nature of deep learning, the proposed method is capable of solving high-dimensional problems on complicated domains with multiple solutions, while existing methods focus on merely one or two-dimensional regular domains and are more expensive in operation counts. Numerical experiments also demonstrate that the proposed method could find more solutions than exiting methods.
\end{abstract}



\pagestyle{myheadings}
\thispagestyle{plain}
\markboth{Structure Probing Neural Network Deflation}{Structure Probing Neural Network Deflation}

\section{Introduction}
\subsection{Problem statement}
Nonlinear differential equations are ubiquitous in various important physical models such as fluid dynamics, plasma physics, solid mechanics, and quantum field theory \cite{Ferris1997,Chyan2000,Davis2001,Kilbas2006}, as well as chemical and biological models \cite{Williams1982,Clark2003}. Solving nonlinear differential equations has been a very challenging problem especially when it is important to find multiple distinct solutions. The nonlinearity of the differential equation may cause traditional iterative solvers to stop at a spurious solution if the initial guess is not close to a physically meaningful solution. When multiple distinct solutions are of interest, a naive strategy is to try different initial guesses as many as possible so that iterative solvers can return distinct solutions as many as possible. However, most of the initial guesses would lead to either spurious solutions or repeated solutions, making this approach usually time-consuming and inefficient unless a priori estimate of the solutions is available.

{ Neural network-based optimization has become a powerful tool for solving nonlinear differential equations, dating back to 1980s \cite{ODE1989} and 1990s \cite{Lee1990,Gobovic1994,Dissanayake1994,Lagaris1998}, and recently revisited in high-dimensional spaces \cite{Han2018,Berg2018,Sirignano2018,Cai2019,Zang2020,Gu2020,RAISSI2019686,Khoo2017SolvingPP,Karpatne2017,Yucesan2020,Pan2020}.} As a form of nonlinear parametrization through compositions of simple functions \cite{IanYoshuaAaron2016}, deep neural networks (DNNs) can efficiently approximate various useful classes of functions without or lessening the curse of dimensionality \cite{barron1993,Hadrien,bandlimit,poggio2017,2019arXiv190609477Y,2020arXiv200103040L,MO,HJKN19_814} and achieve  exponential approximation rates  \cite{yarotsky2017,bandlimit,2020arXiv200103040L,DBLP:journals/corr/LiangS16,DBLP:journals/corr/abs-1807-00297,Opschoor2019}. Therefore, applying DNNs to parametrize the solution space of differential equations (including boundary value problems, initial value problems, and eigenvalue problems) and seeking a solution via energy minimization from variational formulation have become a popular choice, e.g., the least-square method \cite{Berg2018,Sirignano2018,Huang2019} as a special case of variational formulation, the Ritz method \cite{EYu2018}, the Nitsche method \cite{Liao2019}.

However, neural network-based optimization usually can only find the smoothest solution with the fastest decay in the frequency domain due to the implicit regularization of network structures and the stochastic gradient descent (SGD) for solving the minimization problem, no matter how the initial guess is randomly selected. It was shown through the frequency principle of neural networks \cite{DBLP:journals/corr/abs-1905-10264,DBLP:journals/corr/abs-1905-07777,DBLP:journals/corr/abs-1906-09235} and the neural tangent kernel \cite{cao2020towards} that neural networks have an implicit bias towards functions that decay fast in the Fourier domain and the gradient descent method tends to fit a low-frequency function better than a high-frequency function. Through the analysis of the optimization energy landscape of SGD, it was shown that SGD with small batches tends to converge to the flattest minimum \cite{Neyshabur2017,Lei2018,xiaowu}. { Though the above optimization and generalization analysis work only for regression problems, they can be generalized to PDE problems. Recently in \cite{LuoYang2020}, the optimization convergence and generalization analysis of two-layer neural networks for general second-order linear PDEs with variable coefficients on a bounded domain in an arbitrary dimension has been investigated. Global convergence of the gradient descent optimization in the over-parametrization regime is proved using neural tangent kernels and the generalization error with a regularized loss using a Barron-norm is analyzed. Later in \cite{wang2020pinns}, the neural tangent kernel of network-based PDE solvers using two-layer neural networks for one-dimensional Poisson equation is also discussed including the analog of the spectral bias for regression problems proprosed in \cite{cao2020towards}.} Therefore, designing an efficient algorithm for neural network-based optimization to find distinct solutions as many as possible is a challenging problem.

To tackle the challenging problem just above and find distinct solutions as many as possible, we propose a network-based structure probing deflation method in this paper. The key idea of the deflation method is to introduce deflation operators built with known solutions to regularize deep learning optimization, making known solutions no longer local minimizers of the optimization energy landscape while preserving unknown solutions as local minimizers. In particular, we introduce a deflation functional mapping known solutions to infinity. We multiply this deflation functional to the original optimization loss function, then the known solutions will be removed from consideration and unknown solutions can be found by optimizing the regularized loss function via SGD. Furthermore, to facilitate the convergence of SGD, we propose special network structures incorporating boundary conditions of differential equations to simplify the optimization loss function. Finally, a novel structure-probing algorithm is proposed to initialize the deflation optimization making it more powerful to identify distinct solutions with desired structures.

As a general framework, the deflation method can be applied to all neural network-based optimization methods for differential equations. In this paper, we will take the example of boundary value problem (BVP) and the least-square method without loss of generality. The generalization to other problems and methods is similar. Consider the boundary value problem (BVP)
\begin{equation}\label{01}
\begin{split}
&\mathcal{D}u(\bm{x})=f(u(\bm{x}),\bm{x}),\text{~in~}\Omega,\\
&\mathcal{B}u(\bm{x})=g(\bm{x}),\text{~on~}\partial\Omega,
\end{split}
\end{equation}
where $\mathcal{D}: \Omega\rightarrow \Omega$ is a differential operator that is either linear or nonlinear, $f(u(\bm{x}),\bm{x})$ can be a nonlinear function in $u$, $\Omega$ is a bounded domain in $\mathbb{R}^d$, and $\mathcal{B}u=g$ characterizes the boundary condition. Other types of problems like initial value problems can also be formulated as a BVP as discusssed in \cite{Gu2020}. Then least-square method seeks a solution $u(\bm{x};\bm{\theta})$ as a neural network with a parameter set $\bm{\theta}$ via the following optimization problem
\begin{equation}\label{eqn:loss}
\underset{\bm{\theta}}{\min}~L_\text{LS}:=\|\mathcal{D}u(\bm{x};\bm{\theta})-f(u,\bm{x})\|_{L^2(\Omega)}^2+\lambda\|\mathcal{B}u(\bm{x};\bm{\theta})-g(\bm{x})\|_{L^2(\partial\Omega)}^2,
\end{equation}
where $L_\text{LS}$ is the loss function measuring the $L^2$ norms of the residual $\mathcal{D}u(\bm{x};\bm{\theta})-f(u,\bm{x})$ and the boundary residual $\mathcal{B}u(\bm{x};\bm{\theta})-g(\bm{x})$, and $\lambda>0$ is a regularization parameter.

As we shall see, the neural network deflation method enjoys four main advantages compared to traditional deflation methods not based on deep learning:
\begin{itemize}
\item Numerical examples show that the network-based method can identify more solutions than other existing methods, e.g., see Test Case 5 in Section \ref{sec:results}.
\item The network-based method can be applied to solve high-dimensional nonlinear differential equations with multiple solutions while existing methods are only applicable to low-dimensional problems. For example, there is a $6$-dimensional Yamabe's equation in Test Case 6 in Section \ref{sec:results}.
\item The network-based method can be applied to problems with complex domains due to the flexibility of neural network parameterization, e.g., see Test Cases 5 \& 6 in Section \ref{sec:results}.
\item As we shall discuss in Section \ref{sec_complexity}, the network-based method admits lower computational complexity in each iteration compared to existing methods like the original deflation method in \cite{Farrell2015}.
\end{itemize}

\subsection{Related work}
The deflation technique is traced back to the last century for identifying distinct roots of scalar polynomials \cite{Wilkinson1963}. This technique was extended to find roots of systems of nonlinear algebraic equations by Brown and Gearhart in \cite{Brown1971}, where deflation matrices were constructed with old roots to transform the residual of a system of nonlinear algebraic equations so that iterative methods applied to the new residual will only converge to a new root. In \cite{Farrell2015}, Ferrell et al. extended the theoretical framework of Brown and Gearhart \cite{Brown1971} to the case of infinite-dimensional Banach spaces with new classes of deflation operators, enabling the Newton-Kantorovitch iteration to converge to several distinct solutions of nonlinear differential equations even with the same initial guess.

Another well-established method for distinct solutions of differential equations is based on the numerical continuation \cite{AG2003, AG1993, CLP1975, C1979}, where the basic idea of which is to transform the known solutions of a simple start system gradually to the desired solutions of a difficult target system. For example, \cite{Morgan1989} proposed coefficient-parameter polynomial continuation for computing all geometrically isolated solutions to polynomial systems. \cite{Hao2014} put forward a bootstrapping approach for computing multiple solutions of differential equations using a homotopy continuation method with domain decomposition to speed up computation. For more examples of homotopy-based methods and theory in the literature, the reader is referred to \cite{Liao2012}.

The third kind of methods to identify distinct solutions of nonlinear systems is the numerical integration of the Davidenko differential equation associated with the original nonlinear problem \cite{B1972, D1953}. The basic idea is to introduce an artificial time parameter $s$ such that solving the original nonlinear equation $F(u(\bm{x}))=0$ to identify a solution $u_0(\bm{x})$ is equivalent to finding a steady state solution of a time-dependent nonlinear equation $\frac{d F(u(s,\bm{x}))}{ds}+F(u(s,\bm{x}))=0$, which provides a gradient flow of $u(s,\bm{x})$. The gradient flow forms an ordinary differential equation with a solution converging to a solution to the original problem, i.e., $\lim_{s\rightarrow \infty}u(s,\bm{x}) = u_0(\bm{x})$. This method is indeed a broad framework containing the Newton's method as a special example.

\subsection{Organization}
This paper is organized as follows. In Section \ref{sec:RMDE}, we will review the fully connected feed-forward neural network, introduce the formulation of the least-square method for BVP, and design special network structures for four types of boundary conditions. In Section \ref{sec:NND}, the detailed formulation and implementation of the proposed method will be presented. In Section \ref{sec:SP}, the structure probing initialization is introduced. Various numerical experiments are provided in Section \ref{sec:results} to verify the efficiency of the proposed method. Finally, we conclude this paper in Section \ref{sec:conclusion}.

\section{Network-based Methods for Differential Equations}
\label{sec:RMDE}

In this section, we introduce the network-based least-square method based on fully connected feed-forward neural networks and \eqref{eqn:loss} for solving the BVP \eqref{01}. Moreover, special network structures for common boundary conditions are introduced to simplify the loss function in \eqref{eqn:loss} to facilitate the convergence to the desired PDE solution. Vectors are written in bold font to distinguish from scalars in our presentation.

\subsection{Fully connected feed-forward neural network (FNN)}\label{sec:dnn}
FNNs are one of the most popular DNNs and are widely applied to network-based methods for differential equations. Mathematically speaking, for a fixed nonlinear activation function $\sigma$, FNN is the composition of $L$ simple nonlinear functions, called hidden layer functions, in the following formulation:
\[
    \phi(\bm{x};\bm{\theta}):=\bm{a}^T \bm{h}_L \circ \bm{h}_{L-1} \circ \cdots \circ \bm{h}_{1}(\bm{x})\quad \text{for } \bm{x}\in\mathbb{R}^d,
\]
where $\bm{a}\in \mathbb{R}^{N_L}$; $\bm{h}_{\ell}(\bm{x}_{\ell}):=\sigma\left(\bm{W}_\ell \bm{x}_{\ell} + \bm{b}_\ell \right)$ with $\bm{W}_\ell \in \mathbb{R}^{N_{\ell}\times N_{\ell-1}}$ and $\bm{b}_\ell \in \mathbb{R}^{N_\ell}$ for $\ell=1,\dots,L$. With the abuse of notations, $\sigma(\bm{x})$ means that $\sigma$ is applied entry-wise to a vector $\bm{x}$ to obtain another vector of the same size. Usual choices of $\sigma$ include the rectified linear unit (ReLU) function $\sigma(x)=\max\{x,0\}$, its cubic polynomial $\sigma(x)=\max\{x^3,0\}$, a hyperbolic tangent function $\sigma
(x)=\tanh(x)$, etc. $N_\ell$ is the width of the $\ell$-th layer and $L$ is the depth of the FNN. $\bm{\theta}:=\{\bm{a},\,\bm{W}_\ell,\,\bm{b}_\ell:1\leq \ell\leq L\}$ is the set of all parameters in $\phi$ to determine the underlying neural network. Other kinds of neural networks are also suitable in our proposed methods, but we will adopt FNNs for simplicity.

\subsection{Least-square method}
The least-square method is an optimization approach to solve general differential equations. Specifically, let $u(\bm{x}; \bm{\theta})$ be a neural network to approximate the solution $u(\bm{x})$ of BVP \eqref{01}, then the least-square method is formulated as
\begin{equation}\label{03}
\underset{\bm{\theta}}{\min}~L_\text{LS}(\bm{\theta}):=\|\mathcal{D}u(\bm{x};\bm{\theta})-f(\bm{x})\|_{L^2(\Omega)}^2+\lambda\|\mathcal{B}u(\bm{x};\bm{\theta})-g(\bm{x})\|_{L^2(\partial\Omega)}^2,
\end{equation}
where $L_\text{LS}$ is the loss function measuring the weighted magnitude of the differential equation residual $\mathcal{D}u(\bm{x};\bm{\theta})-f(\bm{x})$ and the boundary residual $\mathcal{B}u(\bm{x};\bm{\theta})-g(\bm{x})$ in the sense of $L^2$-norm with a weight parameter $\lambda>0$.

The goal of \eqref{03} is to find an appropriate set of parameters $\theta$ such that the network $u(\bm{x};\bm{\theta})$ minimizes the loss $L_\text{LS}$. If the loss $L_\text{LS}$ is minimized to zero with some $\bm{\theta}$, then $u(\bm{x};\bm{\theta})$ satisfies $\mathcal{D}u(\bm{x};\bm{\theta})-f(\bm{x})=0$ in $\Omega$ and $\mathcal{B}u(\bm{x};\bm{\theta})-g(\bm{x})=0$ on $\partial\Omega$, implying that $u(\bm{x};\bm{\theta})$ is exactly a solution of \eqref{01}. If $L_\text{LS}$ is minimized to a nonzero but small positive number, $u(\bm{x};\bm{\theta})$ is close to the true solution as long as \eqref{01} is well-posed (e.g. the elliptic PDE with Neumann boundary condition, see Theorem 4.1 in \cite{Gu2020}).

In general, the optimization problem \eqref{03} is solved by stochastic gradient descent (SGD) method or its variants (e.g. Adagrad \cite{Duchi2011}, Adam \cite{Kingma2014} and AMSGrad \cite{Reddi2019}) in the deep-learning framework. The optimization and mesh-free setting of the least-square method with neural networks admit several advantageous features that led to its great success and popularity including but not limited to 1) the capacity to solve high-dimensional problems; 2) the flexibility to solve equations of various forms on complicated problem domains; 3) the simple and high-performance implementation with automatic differential programming in existing open-source software.

\subsection{Special network structures for boundary conditions}\label{sec:special_network}
In numerical implementations, the least-square loss function in \eqref{03} relies on the selection of a suitable weight parameter $\lambda$ and a suitable initial guess. If $\lambda$ is not appropriate, it may be difficult to identify a reasonably good minimizer of \eqref{03}. For instance, in the BVP \eqref{01} with $g\equiv0$, if we solve \eqref{03} by SGD with an initial guess $\bm{\theta}^0$ such that $u(\bm{x};\bm{\theta}^0)\approx 0$, SGD might converge to a local minimizer corresponding to a solution neural network close to a constant zero, which is far away from the desired solution, especially when the differential operator $\mathcal{D}$ is highly nonlinear or $\lambda$ is too large. The undesired local minimizer is due to the fact that the boundary residual $\|\mathcal{B}u(\bm{x};\bm{\theta})-g(\bm{x})\|$ overwhelms the equation residual $\|\mathcal{D}u(\bm{x};\bm{\theta})-f(\bm{x})\|$ in the loss function.

The idea just above motivates us to design special networks $u(\bm x;\bm \theta)$ that satisfy the boundary condition $\mathcal{B}u(\bm{x};\bm{\theta})=g(\bm{x})$ automatically and hence we can simplify the least-square loss function from \eqref{03} to
\begin{equation}\label{07}
\underset{\bm \theta}{\min}~L_\text{LS}(\bm{\theta}):=\|\mathcal{D}u(\bm x;\bm{\theta})-f(\bm x)\|_{L^2(\Omega)}^2.
\end{equation}
As we shall see in the numerical section, our numerical tests show that such simplification can help SGD to converge to desired solutions rather than spurious solutions. The design of these special neural networks depends on the type of boundary conditions. We will discuss four common types of boundary conditions by taking one-dimensional problems defined in the domain $\Omega=[a, b]$ as an example. Network structures for more complicated boundary conditions in high-dimensional domains can be constructed similarly. In what follows, denote by $\hat{u}(x;\bm{\theta})$ a generic deep neural network with trainable parameters $\bm{\theta}$. We will augment $\hat{u}(x;\bm{\theta})$ with several specially designed functions to obtain a final network ${u}(x;\bm{\theta})$ that satisfies $\mathcal{B}u({x};\bm{\theta})=g({x})$ automatically.

\vspace{0.25cm}

  \noindent  \textbf{Case 1. Dirichlet boundary condition} $u(a)=a_0$, $u(b)=b_0$.

In this case, we can introduce two special functions $h_1(x)$ and $l_1(x)$ to augment $\hat{u}(x;\bm{\theta})$ to obtain the final network $u(x;\bm{\theta})$:
\begin{equation}\label{11}
u(x;\bm{\theta}) = h_1(x)\hat{u}(x;\bm{\theta})+l_1(x).
\end{equation}
Note $h_1(x)$ and $l_1(x)$ are chosen such that $u(x;\bm{\theta})$ automatically satisfies the Dirichlet $u(a;\bm{\theta})=a_0, \ u(b;\bm{\theta})=b_0$
no matter what $\bm{\theta}$ is. Then $u(x;\bm{\theta})$ is used to approximate the true solution of the differential equation and is trained through \eqref{07}.

For the purpose, $l_1(x)$ is set as a lifting function which satisfies the given Dirichlet boundary condition, i.e. $l_1(a)=a_0$, $l_1(b)=b_0$; $h_1(x)$ is set as a special function which satisfies the homogeneous Dirichlet boundary condition, i.e. $h_1(a)=0$, $h_1(b)=0$. A straightforward choice of $l_1(x)$ is the linear function given by
\begin{equation*}
l_1(x)=(b_0-a_0)(x-a)/(b-a)+a_0.
\end{equation*}
For $h_1(x)$, we can set it as a (possibly fractional) polynomial with roots $a$ and $b$, namely,
\begin{equation*}
h_1(x) = (x-a)^{p_a}(x-b)^{p_b},
\end{equation*}
with $0<p_a,~p_b\leq1$. To obtain an accurate approximation, $p_a$ and $p_b$ should be chosen to be consistent with the orders of $a$ and $b$ of the true solution, hence no singularity will be brought into the network structure. For regular solutions, we take $p_a=p_b=1$; for singular solutions, $p_a$ and $p_b$ would take fractional values. For instance, in the case of a fractional Laplace equation $(-\Delta)^s u=f$ for $0<s<1$ on the domain $\Omega=[-1, 1]$ with boundary conditions $u(\pm 1)=0$, the true solution $u(x)$ has the property that $u(x)=(x-1)^s(x+1)^sv(x)$ with $v(x)$ as a smooth function \cite{Acosta2016,Dyda2017}. Then in the construction of $u(x;\bm{\theta})$, it is reasonable to choose $h_1(x)=(x-1)^s(x-1)^s$ and $l_1(x)=0$.

\vspace{0.25cm}
\noindent \textbf{Case 2. one-sided condition} $u(a)=a_0$, $u'(a)=a_1$.

Similarly to Case 1, the special network is constructed by $u(x;\bm{\theta})=h_2(x)\hat{u}(x;\bm{\theta})+l_2(x)$, where the lifting function $l_2(x)$ is given by
\begin{equation*}
l_2(x)=a_1(x-a)+a_0,
\end{equation*}
and $h_2(x)$ is set as
\begin{equation}
h_2(x) = (x-a)^{p_a},
\end{equation}
with $1<p_a\leq2$. Such $p_a$ guarantees $h_2(x)\hat{u}(x;\bm{\theta})$ and its first derivative both vanish at $x=a$.

\vspace{0.25cm}
\noindent \textbf{Case 3. mixed  boundary condition} $u'(a)=a_0$, $u(b)=b_0$.

In this case, the special network is constructed by $u(x;\bm{\theta}) = \tilde{u}(x; \bm{\theta})+l_3(x)$ with a lifting function $l_3(x)$ chosen as a linear function satisfying the mixed boundary conditions, e.g.,
\begin{equation*}
l_3(x)=a_0x+b_0-a_0b,
\end{equation*}
and $\tilde{u}(x;\bm{\theta})$ satisfying the homogeneous mixed boundary conditions. In the construction of $\tilde{u}(x;\bm{\theta})$, it is inappropriate to naively take $\tilde{u}(x;\bm{\theta})=(x-a)^{p_a}(x-b)^{p_b}$ with $1<p_a\leq2$ and $0<p_b\leq1$, following the approaches in the preceding two cases, because such $\tilde{u}(x;\bm{\theta})$ satisfies a redundant condition $\tilde{u}(a;\bm{\theta})=0$. Instead, we assume
\begin{equation}\label{04}
\tilde{u}(x;\bm{\theta})=(x-a)^{p_a}\hat{u}(x;\bm{\theta})+c,
\end{equation}
where $1<p_a\leq 2$ and $c$ is a network-related constant to be determined. Clearly, \eqref{04} implies $\tilde{u}'(a;\bm{\theta})=0$, whereas $\tilde{u}(a;\bm{\theta})$ has not been specified. Next, the constraint $\tilde{u}(b;\bm{\theta})=0$ gives $c=-(b-a)^{p_a}\hat{u}(b;\bm{\theta})$. Therefore, the special network for mixed boundary conditions is constructed via
\begin{equation}\label{12}
u(x;\bm{\theta}) = (x-a)^{p_a}\hat{u}(x;\bm{\theta})-(b-a)^{p_a}\hat{u}(b;\bm{\theta})+l_3(x).
\end{equation}

\vspace{0.25cm}
   \noindent \textbf{Case 4. Neumann boundary condition} $u'(a)=a_0$, $u'(b)=b_0$.

Similarly to Case 3, we construct the network by $u(x;\bm{\theta}) = \tilde{u}(x; \bm{\theta})+l_4(x)$ with a lifting function $l_4(x)$ satisfying the Neumann boundary condition given by
\begin{equation*}
l_4(x)=\frac{(b_0-a_0)}{2(b-a)}(x-a)^2+a_0x.
\end{equation*}
And $\tilde{u}(x;\bm{\theta})$ satisfying the homogeneous Neumann boundary condition is assumed to be
\begin{equation}\label{05}
\tilde{u}(x;\bm{\theta})=(x-a)^{p_a}\check{u}(x;\check{\bm{\theta}})+c_1,
\end{equation}
where $1<p_a\leq2$, $\check{u}(x;\check{\bm{\theta}})$ is an intermediate network to be determined later, and $c_1$ is a network parameter to be trained together with $\check{\bm{\theta}}$. It is easy to check that $\tilde{u}'(a;\bm{\theta})=0$. Next, by the constraint $\tilde{u}'(b;\bm{\theta})=p_a(b-a)^{p_a-1}\check{u}(b;\check{\bm{\theta}})+(b-a)^{p_a}\check{u}'(b;\check{\bm{\theta}})=0$, we have
\begin{equation*}
p_a\check{u}(b;\check{\bm{\theta}})+(b-a)\check{u}'(b;\check{\bm{\theta}}) =0,
\end{equation*}
which can be reformulated as
\begin{equation*}
\Big(\exp(\frac{p_a x}{b-a})\check{u}(x;\check{\bm{\theta}})\Big)'_{x=b}=0.
\end{equation*}
Therefore, we have
\begin{equation}\label{06}
\exp(\frac{p_a x}{b-a})\check{u}(x;\check{\bm{\theta}})=(x-b)^{p_b}\hat{u}(x;\hat{\bm{\theta}})+c_2,
\end{equation}
where $1<p_b\leq 2$ and $c_2$ is another network parameter to be trained together with $\hat{\bm{\theta}}$. Finally, by combining \eqref{05} and \eqref{06}, we obtain the following special network satisfying the given Neumann condition, i.e.
\begin{equation}\label{13}
u(x;\bm{\theta})=\exp(\frac{p_a x}{a-b})(x-a)^{p_a}\big((x-b)^{p_b}\hat{u}(x;\hat{\bm{\theta}})+c_2\big)+c_1+l_4(x),
\end{equation}
where $\bm{\theta}=\{\hat{\bm{\theta}},c_1,c_2\}$.

{ Finally, we would like to remark that it is difficult to construct special neural networks to automatically satisfy boundary conditions when the PDE domain is irregular. In this case, the conventional penalty method in \eqref{03} is more preferable. Though we will show in our numerical experiments that special neural networks satisfying boundary conditions are better than penalty methods to identify distinct solutions. This does not exclude the possibility that penalty methods, or other advanced optimization algorithms for constrained optimization, can also work well with well-tuned parameters.}

\section{Neural Network Deflation}
\label{sec:NND}

In this section, we propose the general formulation, the detailed implementation, and the computational complexity of the proposed method. As we shall see, our method is easy to implement on high-dimensional and complex domains with a lower computational cost per iteration than other traditional deflation methods.

\subsection{Formulation}
A nonlinear BVP \eqref{01} might have multiple distinct solutions and each solution is a minimizer of the corresponding network-based optimization, say
\begin{equation}\label{02}
\underset{\bm \theta}{\min}~L(u(\bm x;\bm \theta)),
\end{equation}
where $L$ is a generic loss function for solving differential equations. One example of $L$ is the residual loss in \eqref{07}, and $L$ can also be other loss functions. Due to the implicit regularization of SGD and neural networks, only local minimizers in flat energy basins are likely to be found. Hence, no matter how to initialize the SGD and how to choose hyper-parameters, usually, only a few solutions can be found by minimizing \eqref{02} directly.

{ The neural network deflation is therefore introduced, the main idea of which is to construct a modified loss function $L_\text{ND}$ with two properties: First, a candidate minimizer of $L_\text{ND}$ is also a minimizer of $L$. Second, the minimizers that are already found by the network-based optimization \eqref{02} will not be minimizers of $L_\text{ND}$ again. Following this idea, $L_\text{ND}$ is constructed by multiplying $L$ with a deflation term introduced in \cite{Farrell2015} such that the energy landscape of $L$ is modified. Specifically, suppose the minimum value of $L$ is zero. Let $u_k(\bm x)$ ($k=1, \cdots, K$) be the solutions already found by \eqref{02}, then the neural network deflation is formulated as the following optimization problem,
\begin{equation}\label{09}
\underset{\bm \theta}{\min}~L_\text{ND}:=\Big(\overset{K}{\underset{k=1}{\sum}}\|u(\bm x;\bm \theta)-u_k(\bm x)\|_{L^2(\Omega)}^{-p_k}+\alpha\Big) L(u(\bm x;\bm \theta)),
\end{equation}
where $p_k$ are positive powers for $k=1,\cdots,K$, and $\alpha>0$ is a shift constant. Here, we name $u_k(\bm x)$ ($k=1, \cdots, K$) as deflation sources. Indeed, the modified loss function $L_\text{ND}$ satisfies the two properties discussed above. First, any minimizer of $L_\text{ND}$ such that $L_\text{ND}=0$ also ensure $L=0$ and, hence, is also a minimizer of $L$. Second, for all $k=1, \cdots, K$, the term $\|u(\bm x;\bm \theta)-u_k(\bm x)\|_{L^2(\Omega)}^{-p_k}$ acts a penalty term that excludes $u_k$ as a minimizer, since it approaches infinity as $u$ goes to $u_k$. The introduction of a positive $\alpha$ is help to eliminate spurious solutions in practice. If $\alpha=0$, the modified loss function $L_\text{ND}$ would approach zero when $u$ is far from all $u_k$'s, which leads to many spurious solutions. For a more detailed discussion of the deflation term, the reader can refer to \cite{Farrell2015}. }

\subsection{Deflation with a varying shift} The original deflation operator introduced in \cite{Farrell2015} fixes the shift $\alpha$ in \eqref{09} as a constant. In this paper, we propose a new variant of deflation operators with a varying shift $\alpha$ along with the SGD iteration. Note that when $\alpha$ is equal or close to $0$, the deflation term $\sum_{k=1}^{K}\|u(\bm x;\bm \theta)-u_k(\bm x)\|_{L^2(\Omega)}^{-p_k}$ dominates the loss and hence gradient descent tends to converge to what is far away from the known solutions. When $\alpha$ is moderately large, the original loss function $L(u(\bm x;\bm \theta))$ dominates the loss and the gradient descent process tends to converge to a solution with a smaller residual. Therefore, $\alpha$ in this paper is set to be a monotonically increasing function of the SGD iteration. In the early stage, $\alpha$ is chosen to be close to $0$ such that the current solution will be pushed away from known solutions. During this stage, a large learning rate is preferable. In the latter stage when the current solution is roughly stable, $\alpha$ is set to be large and a small learning rate is used to obtain a small residual loss.

{ In practice, one heuristic choice is to increase $\alpha$ exponentially with a linearly growing power when the iteration increases. For example, in the $n$-th iteration, $\alpha$ is set as $\alpha_n$ defined below
\begin{equation}\label{15}
\alpha_n=10^{p_0+n(p_1-p_0)/N_\text{I}},
\end{equation}
where $p_0$ and $p_1$ are two prescribed powers with $p_0\leq p_1$, and $N_\text{I}$ is the total number of iterations. Note that the exponentially varying formula is also widely used in setting the learning rates of SGD.}

\subsection{Discretization}
In the implementation, the continuous loss functions in \eqref{07} and \eqref{09} are approximately evaluated by stochastic sampling. The $L^2$-norm can be interpreted as an expectation {of a function of a random} variable $\bm{x}$ in a certain domain. Hence, the expectation is approximated by sampling $\bm{x}$ several times and computing the average function value as an approximant. Let us take $\|u(\bm{x})\|_{L^2(\Omega)}$ as an example. We generate $N_\text{p}$ random samples $\bm x_i$, $i=1,\cdots,N_\text{p}$, which are uniformly distributed in $\Omega$. Denote $\bm X:=\{\bm x_i\}_{i=1}^{N_\text{p}}$, then $\|u(\bm{x})\|_{L^2(\Omega)}$ is evaluated as the discrete $L^2$-norm denoted as $\|u(\bm{x})\|_{L^2(\bm X)}$ via
\begin{equation}\label{14}
\|u(\bm x)\|_{L^2(\bm X)}:=\Big(\frac{1}{N_\text{p}}\underset{\bm x_i\in \bm X}{\sum}|u(\bm x_i)|^2\Big)^{\frac{1}{2}}.
\end{equation}

The discretization technique above is applied to discretize the $L^2$-norms in all loss functions in this paper. In the $n$-th iteration of gradient descent, assuming that the shift $\alpha$ is set to be $\alpha_n$ and the set of random samples is denoted as $\bm X_n$, the discrete deflation loss function is calculated by
\begin{equation*}
\widehat{L}_\text{ND}^{(n)}(\bm{\theta}):=\Big(\overset{K}{\underset{k=1}{\sum}}\|u(\bm x;\bm \theta)-u_k(\bm x)\|_{L^2(\bm X_n)}^{-p_k}+\alpha_n\Big) \widehat{L}(u(\bm x;\bm \theta)),
\end{equation*}
where $\widehat{L}(u(\bm x;\bm \theta))$ is a discrete approximation to $L(u(\bm x;\bm \theta))$ using the same set of samples, e.g.,
\[
\widehat{L}(u(\bm x;\bm \theta))=\|\mathcal{D}u(\bm x;\bm{\theta})-f(\bm x)\|_{L^2(\bm{X}_n)}^2
\]
when the least-square loss in \eqref{07} is applied.
Then the network parameter $\bm \theta$ is updated by
\begin{equation*}
\bm \theta\leftarrow\bm \theta-\tau_n\nabla_{\bm \theta}\widehat{L}_\text{ND}^{(n)}(\bm{\theta}),
\end{equation*}
where $\tau_n>0$ is the learning rate in the $n$-th iteration. In our implementation, $\bm{X}_n$ is renewed in every iteration. { Note that the gradient of the loss function can be evaluated using PyTorch built-in function autograd that essentially compute{s} the gradient using a sequence of chain rules, since the network is the composition of several simple functions with explicit formulas.}

\subsection{Computational Complexity}\label{sec_complexity}
Let us estimate the computational complexity of the SGD algorithm for deflation optimization \eqref{09} with least-square loss function \eqref{07}. Recall that $N_\text{p}$ denotes the number of random samples in each iteration. Assume that the FNN has $L$ layers and $N$ neurons in each hidden layer. Note that evaluating the FNN or computing its derivative with respect to its parameters or input $\bm{x}$ via the forward or backward propagation takes $O(d N+L N^2)$ FLOPS (floating point operations per second) for each sample $\bm{x}$. Moreover, as in most existing approaches, we assume $f(\bm x)$ in the BVP can be evaluated with $O(d)$ FLOPS for a single $\bm x$. Therefore, $L(u(\bm x;\bm \theta))$ in \eqref{07} and its derivative $\nabla_{\bm\theta} L(u(\bm x;\bm \theta))$ can be calculated with $O(N_\text{p}(d N+L N^2))$ FLOPS using the discrete $L^2$-norm in \eqref{14}, if the differential operator $\mathcal{D}$ is evaluated through finite difference approximation. Similarly, assuming the number of known solutions $K$ is $O(1)$ and the known solutions $\{u_k(\bm x)\}_{k=1}^K$ are stored as neural networks of width $N$ and depth $L$, then the deflation factor and its gradient with respect to ${\bm\theta}$ can also be calculated with $O(N_\text{p}(d N+L N^2))$ FLOPS. Finally, the total complexity in each gradient descent iteration of the deflation optimization is $O(N_\text{p}(d N+L N^2))$.

In existing methods \cite{Farrell2015,Croci2015,Adler2016}, a given nonlinear differential equation is discretized via traditional discretization techniques, e.g. FDM and FEM, resulting in a nonlinear system of algebraic equations. The solutions of the system of algebraic equations provide numerical solutions to the original nonlinear differential equation. By multiplying different deflation terms to the nonlinear system of algebraic equations, existing methods can identify distinct solutions via solving the deflated system by Newton's iteration. The number of algebraic equations $N_e$ derived by FDM is exactly the number of grid points; and the number of equations derived by FEM is exactly the number of trial functions in the Galerkin formulation.

Now we compare neural network deflation with existing deflation methods in \cite{Farrell2015,Croci2015,Adler2016} in terms of the computational complexity under the assumption that the degrees of freedom of these methods are equal, i.e., the number of grid points or trial functions in existing methods is equal to the number of parameters in the neural network deflation, which guarantees that these methods have almost the same accuracy to find a solution. Denote the degree of freedom of these methods by $W$. Then by the above discussion, we have $W=N_e=O(d N + L N^2)$. Therefore, the total computational complexity in each iteration is $O(N_p W)$, where $N_p$ is usually chosen as a hyper-parameter much smaller than $W$. In existing methods, the Jacobian matrix in each Newton's iteration is a low-rank matrix plus a sparse matrix of size $W$ by $W$. Typically, each iteration of Newton's method requires solving a linear system of the Jacobian matrix, which usually requires $O(W^2)$ FLOPS. If a good preconditioner exists or a sparse direct solver for inverting the Jacobian matrix exists, the operation count may be reduced. Consequently, the total complexity in each iteration of existing methods would be more expensive than the neural network deflation depending on the performance of preconditioners.

\section{Structure Probing Initialization}
\label{sec:SP}

The initialization of parameters plays a critical role in training neural networks and has a significant impact on the ultimate performance. In the training of a general FNN, network parameters are usually randomly initialized using normal distributions with zero-mean. One popular technique is the Xavier initialization \cite{Glorot2010}: for each layer $\ell$, the weight matrix $\bm{W}_\ell\in\mathbb{R}^{N_\ell\times N_{\ell-1}}$ is chosen randomly from a normal distribution with mean $0$ and variance $1/N_{l-1}$; the bias vector $\bm{b}_\ell$ is initialized to be zero. As a variant of Xavier initialization, the He initialization \cite{He2015} takes a slightly different variance $2/(N_{l-1}+N_l)$ for $W_\ell$ and $2/N_{l-1}$ for $b_\ell$. In general, FNNs initialized randomly have a smooth function configuration, and hence their Fourier transform coefficients decay quickly.

The least-squares optimization problem, either for regression problems or solving linear partial differential equations,  with over-parameterized FNNs and random initialization admits global convergence by gradient descent with a linear convergence rate \cite{DBLP:journals/corr/abs-1806-07572,Du2018,Zhu2019,Chen1,LuoYang2020}. However, the speed of convergence depends on the spectrum of the target function. The training of a randomly initialized DNN tends to first capture the low-frequency components of a target solution quickly. The high-frequency fitting error cannot be improved significantly until the low-frequency error has been eliminated, which is referred to as F-principle \cite{Xu2019}.
Related works on the learning behavior of DNNs in the frequency domain is further investigated in \cite{DBLP:journals/corr/abs-1906-09235,DBLP:journals/corr/abs-1905-07777,DBLP:journals/corr/abs-1905-10264,cao2020towards}. In the case of nonlinear differential equations where multiple solutions exist, these theoretical works imply that deep learning-based solvers converge to solutions in the low-frequency domain unless the DNN is initialized near a solution with high-frequency components.

The discussion just above motivates us to propose the structure probing initialization that helps the training converge to multiple structured solutions. The structure probing initialization incorporates desired structures in the initialization and training of DNNs. For example, to obtain oscillatory solutions of a differential equation, we initialize the DNN with high-frequency components to make the initialization closer to the desired oscillatory solution. During the optimization process, the magnitudes of these high-frequency components will be optimized to fit the desired solution. One choice to probe an oscillatory solution is to take a linear combination of structure probing functions with various frequencies, e.g., $\{\xi_j(\bm x)=e^{\text{i}\bm{k}_j \cdot \bm x}, |\bm{k}_j|=j, j=1, \cdots, J\}$ with $\bm{k}_j$ randomly selected. Then the following network $u_J$ with a set of random parameters $\bm \theta$ can serve as an oscillatory initial guess:
\begin{equation}\label{uj}
u_J(\bm x;\bm \theta_J)=u(\bm x;\bm \theta)+\underset{j=1}{\overset{J}{\sum}}{c}_j\xi_j(\bm x),
\end{equation}
where $\bm{\theta}_J:=\{\bm\theta,\{c_j\}_{j=1}^J\}$ is trainable after initialization. In the initialization, $\{c_j\}$ are set as random numbers or manually determined hyper-parameters with large magnitudes. The idea of adding planewaves has been applied in \cite{cai2019phase,cai2019phasednn} to obtain high-frequency solutions. But the goal and detailed formulations are different. Instead of planewaves, radial basis functions are also a popular structure in the solution of differential equations. In this case, we can choose $\{\xi_j(\bm x)=\sin(j\pi|\bm x|), j=1, \cdots, J\}$ for example. The idea of structure probing initialization is not limited to the above two types of structures and is application dependent.

The above paragraph has sketched out the main idea of the structure probing initialization. Now we are ready to discuss its special cases when we need to make the structure probing network $u_J$ in \eqref{uj} satisfy the boundary condition $\mathcal{B}u_J=g$ in the BVP \eqref{01}, which is important for the convergence of deep learning-based solvers as discussed in Section \ref{sec:special_network}. For this purpose, we first construct a special network $u(\bm x;\bm \theta)$ such that $\mathcal{B}u(\bm x;\bm\theta)=g$ by the approaches described in Section \ref{sec:special_network}. Next, the structured probing functions $\{\xi_j(\bm x)\}$ are specifically chosen to satisfy $\mathcal{B}\xi_j(\bm x)=0$ for each $j$. As an example, let us take the one-dimensional mixed boundary condition on $[a,b]$:
\begin{equation}
u'(a)=a_0,\quad u(b)=b_0
\end{equation}
for any constants $a_0$ and $b_0$. Then a feasible choice of $\xi_j(\bm x)$ is $\xi_j(x)=\cos(\frac{(2j-1)\pi (x-a)}{2(b-a)})$. Finally, it is easy to check that $\mathcal{B}u_J(\bm x;\bm \theta)=g$.

\section{Numerical Examples}\label{sec:results}
In this section, several numerical examples are provided to show the performance of network-based structure probing deflation in solving BVP \eqref{01}. We choose the least-square loss function as the general loss function $L(u(\bm x;\bm \theta))$ in \eqref{09}, then the neural network deflation is formulated as
\begin{equation}\label{08}
\underset{\bm \theta}{\min}~L_\text{ND}(\bm{\theta}):=\Big(\overset{K}{\underset{k=1}{\sum}}\|u(\bm x;\bm \theta)-u_k(\bm x)\|_{L^2(\Omega)}^{-p_k}+\alpha\Big) \|\mathcal{D}u(\bm x;\bm \theta)-f(\bm x)\|_{L^2(\Omega)}^2,
\end{equation}
where $u(\bm x;\bm \theta)$ is the neural network of the approximate solution to be determined. Remark that the optimization problem can also be formulated by other optimization-based methods instead of least squares.

To verify the effectiveness of special networks that satisfy boundary conditions automatically, we use the deflation without the special network for boundary conditions as a comparison, where the loss function of the deflation becomes
\begin{multline}\label{16}
\underset{\bm \theta}{\min}~L_\text{ND}(\bm{\theta}):=\Big(\overset{K}{\underset{k=1}{\sum}}\|u(\bm x;\bm \theta)-u_k(\bm x)\|_{L^2(\Omega)}^{-p_k}+\alpha\Big)\cdot\\ \left(\|\mathcal{D}u(\bm{x};\bm{\theta})-f(u,\bm{x})\|_{L^2(\Omega)}^2+\lambda\|\mathcal{B}u(\bm{x};\bm{\theta})-g(\bm{x})\|_{L^2(\partial\Omega)}^2\right).
\end{multline}

The overall setting for all examples is summarized as follows.
\begin{itemize}
  \item \textbf{Environment.}
  The experiments are performed in Python 3.7 environment. We utilize the PyTorch library for neural network implementation and CUDA 10.0 toolkit for GPU-based parallel computing. One-dimensional examples (Test Case 1-4) are implemented on a laptop and high-dimensional examples (Test Case 5-7) are implemented on a scientific workstation.
  \item \textbf{Optimizer.}
  In all examples, the optimization problems are solved by {\em adam} subroutine from PyTorch library with default hyper parameters. This subroutine implements the Adam algorithm in \cite{Kingma2014}.
  \item \textbf{Learning rate.}
  In each example, the learning rate is set to decay exponentially with linearly decreasing powers. Specifically, the learning rate in the $n$-th iteration is set as
    \begin{equation*}
    \tau_n=10^{q_0+n(q_1-q_0)/N_\text{I}},
    \end{equation*}
    where $q_0>q_1$ are the initial and final powers, respectively, and $N_\text{I}$ denotes the total number of iterations.
  \item \textbf{Network setting.}
  In each example, we construct a special network that satisfies the given boundary condition as discussed in Section \ref{sec:special_network}. The special network involves a generic FNN, denoted by $\hat{u}$. In all examples, we set the depth and width of $\hat{u}$ as fixed numbers $L=3$ and $N=100$. Unless specified particularly, all weights and biases of $\hat{u}$ are initialized by ${\bm W}_l, {\bm b}_l\sim U(-\sqrt{N_{l-1}},\sqrt{N_{l-1}})$. The activation function of $\hat{u}$ is chosen as $\sigma(x):=\max(0,x^3)$.
  \item \textbf{Varying shifts in deflation operators.}
  In one-dimensional examples (Test Case 1-4), using constant shifts is sufficient to find all solutions. In high-dimensional examples (Test Case 5-7), varying shifts will help to find more distinct solutions. In these examples, we set varying shifts according to \eqref{15}.
\end{itemize}

We also summarize the numerical examples in this section in Table \ref{tab: summary} below, which could help the reader to better understand how the extensive numerical examples demonstrate the advantages of our new ideas in this paper: 1) neural network deflation (ND); 2) structure probing initialization (SP); 3) special network for boundary conditions (BC); 4) varying shifts in deflation operators (VS).

\begin{table}[H]
    \small
    \centering
        \begin{tabular}{|c||cccc||c|}
        \hline
          Test Case  & ND & SP & BC & VS  & Justified Ideas   \\
        \hline
          Case 1 & 1/0 & 0 & 1/0 & 0  & ND and BC  \\
          Case 2 & 1/0 & 0 & 1/0 & 0  & ND and BC   \\
          Case 3 & 1/0 & 0 & 1/0 & 0  & ND and BC   \\
          Case 4 & 1/0 & 1/0 & 1/0 & 0  & ND, SP and BC   \\
          Case 5 & 1/0 & 1/0 & 1 & 1/0  & ND, SP, and VS   \\
          Case 6 & 1/0 & 1/0 & 1 & 1/0  & ND, SP, and VS   \\
          Case 7 & 1/0 & 0 & 1 & 1  & ND   \\
        \hline
\end{tabular}
    \caption{Summary of numerical examples and goals. In this table, ``1" represents an idea is used and ``0" means the idea is not used. ``1/0" indicates that a comparison with/without the idea is tested.}
    \label{tab: summary}
\end{table}

In each example, necessary parameters to obtain each solution are listed in a table right next to the example. In these tables, we use $N_\text{p}$, $N_\text{I}$, and $I_\text{lr}$ to denote { the batch size, the number of iterations, and the range of learning rates (i.e. $[10^{q_1},10^{q_0}]$),} respectively. In each iteration of the optimization, $N_\text{p}$ random samples will be renewed. The value of the shift $\alpha$ for each solution found by the deflation is listed in the table as a constant for a fixed $\alpha$ or an interval $[10^{p_0},10^{p_1}]$ for a varying $\alpha$.

\subsection{Numerical tests in one-dimension} First of all, we will provide four numerical tests for problems in one-dimension. These numerical tests show that the proposed neural network deflation works as well as existing methods \cite{Farrell2015,Hao2014}.

\vspace{0.25cm}
\noindent\large{\textbf {Test Case 1.}}
We consider second-order the Painlev{\'e} equation \cite{HOLMES1984,Fornberg2011,Noonburg1995} that seeks $u(x)$ satisfying
\begin{eqnarray}\label{1-D_Painleve}
&&\frac{\td^2u}{\td x^2}=100u^2-1000x, \qquad \text{in}\  \Omega=(0,1),\\
 && u(0)=0,\quad u(1)=\sqrt{10}.\label{bou}
\end{eqnarray}
It has been shown in \cite{Hastings1989} that the Painlev{\'e} equation \eqref{1-D_Painleve}-\eqref{bou} has exactly two solutions, denoted by $u_1$ and $u_2$, which satisfy $u_1'(0)>0$ and $u_2'(0)<0$, respectively.

\begin{figure}
\centering
\includegraphics[scale=0.8]{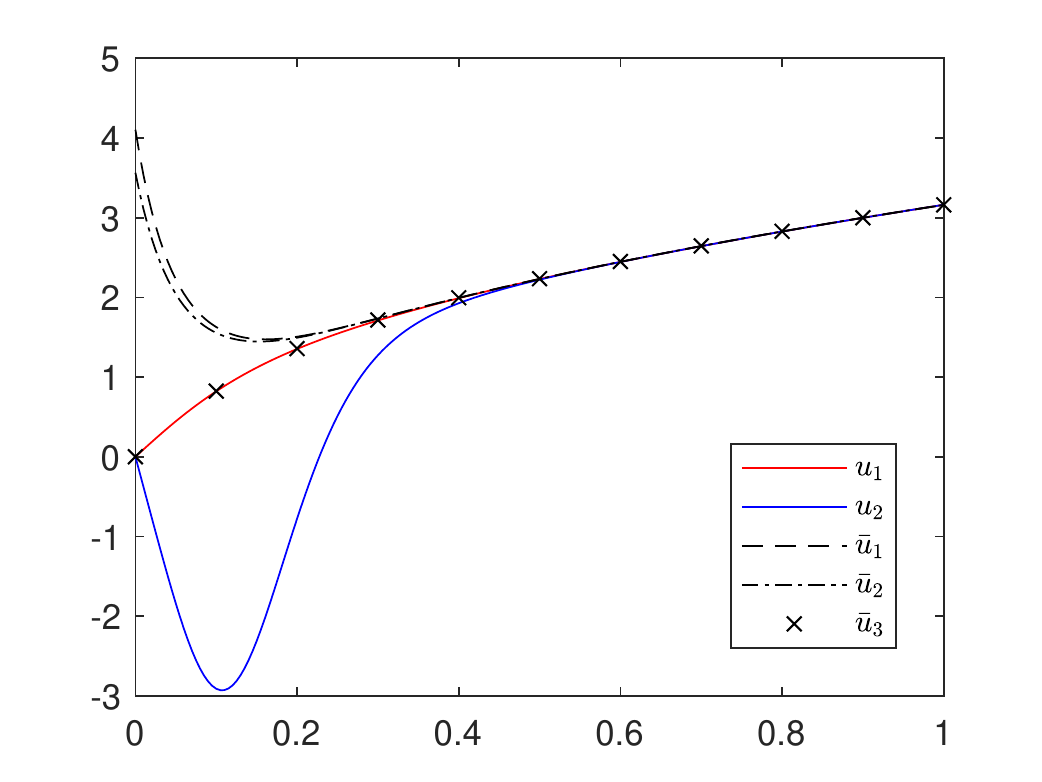}
\caption{\em Identified solutions of the 1-D Painlev{\'e} equations \eqref{1-D_Painleve}-\eqref{bou} by the least squares method and neural network deflation. All correct solutions, $u_1$ and $u_2$, are identified with special networks for boundary conditions. Spurious solutions, $\bar{u}_1$ and $\bar{u}_2$, are found if the special networks are not used. Another solution, $\bar{u}_3$, is found when the deflation fails with an inappropriate power $p_1$.}
\label{Painleve_solutions}
\end{figure}

In our experiments, we take the following special network
\begin{equation}\label{18}
u(x;\bm\theta)=x(x-1)\hat{u}(x;\bm\theta)+\sqrt{10}x
\end{equation}
that automatically satisfies the boundary conditions and use parameters in Table \ref{Painleve_parameters}. The initial guess of $\bm{\theta}$ is randomly initialized as mentioned previously. The first solution $u_1$ is easily found by the least-square method using \eqref{07}, and the second solution $u_2$ is found by deflation with $u_1$ as the deflation source and $p_1=2$. Other parameters associated with these solutions are listed in Table \ref{Painleve_parameters}. Figure \ref{Painleve_solutions} visualizes the identified solutions $u_1$ and $u_2$ with the same function configurations as in \cite{Farrell2015}.

To verify the effectiveness of special networks that satisfy the boundary conditions \eqref{bou}, we use the deflation without special networks for boundary conditions as a comparison. Hence, the loss function is given by \eqref{16} with a solution network $u(\bm{x};\bm{\theta})$ as a generic FNN of the same structure as $\hat{u}$ in \eqref{18}. To show that the results of \eqref{16} are quite independent of the weight $\lambda$, $\lambda=1$ and $\lambda=100$ are used and the corresponding solutions are denoted as $\bar{u}_1$ and $\bar{u}_2$, respectively. As listed in Table \ref{Painleve_parameters}, other parameters to identify these two solutions are the same as those for identifying $u_2$ for a fair comparison. It is clear that these two solutions do not satisfy the boundary condition at the endpoint $x=0$ (see Figure \ref{Painleve_solutions}). This verifies the importance of using special networks that satisfy the boundary conditions automatically.

{ Moreover, we test the effectiveness of the deflation when smaller powers of deflation sources are used. We basically repeat the same experiment as the previous one for computing $u_2$ in Figure \ref{Painleve_solutions}. The only difference is that we use a larger power $p_1=2$ in the previous experiment, but now we use a smaller power $p_1=1$. The solution by nerual network deflation with $p_1=1$ is denoted as $\bar{u}_3$ and visualized in Figure \ref{Painleve_solutions}. Note that $\bar{u}_3$ is almost the same as the deflation source $u_1$ by visual inspection. This result is not surprising even if we have used $u_1$ as the deflation source. The loss function $L_\text{ND}$ in \eqref{09} can still be very small at $u=\bar{u}_3$ even if the deflation term is large, since the loss function $L$ in \eqref{09} can be much smaller than one over the deflation term at $u=\bar{u}_3$ close to $u_1$. This example indicates that an appropriate power $p_1$ is necessary to exclude spurious solutions close to $u_1$. 
}

\begin{table}
\centering
\begin{tabular}{|c|c|c|c|c|}
  \hline
  ~ & $u_1$ & $u_2$ & $\bar{u}_1$ & $\bar{u}_2$\\\hline
  $N_\text{I}$ & 10000 & 10000 & 10000 & 10000\\\hline
  $N_\text{p}$ & 1000 & 1000 & 1000 & 1000\\\hline
  $I_\text{lr}$ & $[10^{-3},10^{-2}]$ & $[10^{-3},10^{-2}]$ & $[10^{-3},10^{-2}]$ & $[10^{-3},10^{-2}]$\\\hline
  deflation source & / & $u_1$ ($p_1=2$) & $u_1$ ($p_1=2$) & $u_1$ ($p_1=2$)\\\hline
  $\alpha$ & / & $1$ & $1$ & $1$ \\
  \hline
\end{tabular}
\caption{\em Parameters for 1-D Painlev{\'e} equations \eqref{1-D_Painleve}-\eqref{bou}. ``/" means the corresponding item is not used (the same as below)}
\label{Painleve_parameters}
\end{table}

\vspace{0.25cm}
\noindent\large{\textbf{Test Case 2.}}
We consider a fourth-order nonlinear BVP that seeks $u$ such that
\begin{eqnarray}\label{1-D_Graef}
&& \frac{\td^4u}{\td x^4}=\beta x(1+u^2) \qquad \text{in}\  \Omega=(0,1),\\
&& u(0)=u'(1)=u''(1)=0,\quad u''(0)-u''(\gamma)=0,\label{bc1}
\end{eqnarray}
where $\beta$ and $\gamma$ are two given constants. Graef et al. \cite{Graef2003_1,Graef2003_2} have proven that the problem \eqref{1-D_Graef}-\eqref{bc1} has at least two positive solutions when $\beta=10$ and $\gamma=1/5$.

The three-point boundary condition \eqref{bc1} is more complicated than usual. Accordingly, we construct the following special network for it,
\begin{equation}\label{21}
u(x;\bm\theta)=(x-1)^3\hat{u}(x;\bm\theta)+\hat{u}(0;\bm\theta)+c_\gamma x(x-1)^3,
\end{equation}
where
\begin{equation}\label{22}
c_\gamma=\frac{1}{-12\gamma^2+18\gamma}\Big(\frac{\td^2}{\td x^2}\big((x-1)^3\hat{u}(x;\bm\theta)\big)|_{x=\gamma}-\frac{\td^2}{\td x^2}\big((x-1)^3\hat{u}(x;\bm\theta)\big)|_{x=0}\Big).
\end{equation}
It can be verified that \eqref{21} indeed satisfies the boundary condition \eqref{bc1} independent of ${\bm \theta}$.

\begin{figure}
\centering
\includegraphics[scale=0.6]{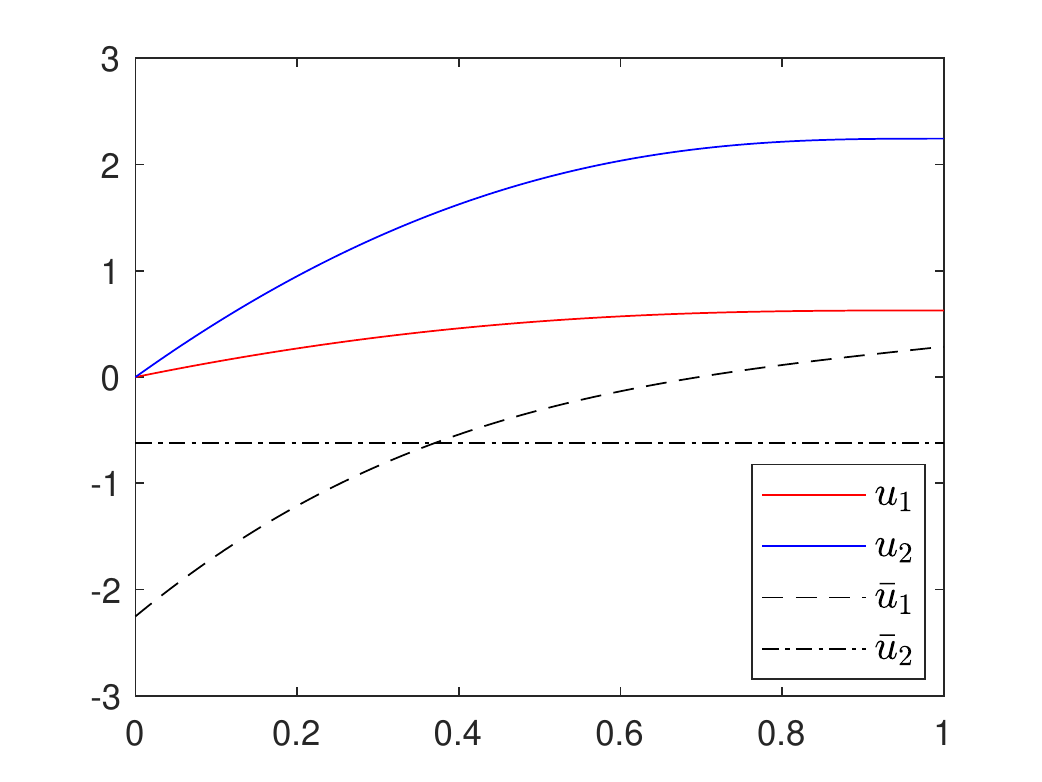}
\caption{\em Identified solutions of the equation \eqref{1-D_Graef}-\eqref{bc1} by least square or neural network deflation. All correct solutions, $u_1$ and $u_2$, are identified with special networks for boundary conditions. Spurious solutions, $\bar{u}_1$ and $\bar{u}_2$, are found if the special networks are not used.}
\label{Graef_solutions}
\end{figure}

In our experiment, we find the first solution, denoted by $u_1$, by applying the least-square method \eqref{07}. With deflation source $u_1$ ($p_1=1$), we find the second solution, denoted by $u_2$, by using the deflation \eqref{08}. The parameters and solutions are demonstrated in Table \ref{Graef_parameters} and Figure \ref{Graef_solutions}.

Similarly, we test the deflation without special networks for boundary conditions as a comparison under the same setting as Test Case 1. We find two solutions, denoted by $\bar{u}_1$ and $\bar{u}_2$, from $\lambda=1$ and $\lambda=100$, respectively (see Figure \ref{Graef_solutions}). It is clear that both solutions are spurious since their configurations do not take the prescribed boundary value $0$ at $x=0$ (see Figure \ref{Graef_solutions}), which implies the effectiveness of using special networks for boundary conditions.

\begin{table}
\centering
\begin{tabular}{|c|c|c|c|c|}
  \hline
  ~ & $u_1$ & $u_2$ & $\bar{u}_1$ & $\bar{u}_2$\\\hline
  $N_\text{I}$ & 5000 & 5000 & 5000 & 5000\\\hline
  $N_\text{p}$ & 1000 & 1000 & 1000 & 1000\\\hline
  $I_\text{lr}$ & $[10^{-3},10^{-2}]$ & $[10^{-2},10^{-1}]$ & $[10^{-2},10^{-1}]$ & $[10^{-2},10^{-1}]$\\\hline
  deflation source & / & $u_1$ ($p_1=1$) & $u_1$ ($p_1=1$) & $u_1$ ($p_1=1$)\\\hline
  $\alpha$ & / & $1$ & $1$ & $1$ \\
  \hline
\end{tabular}
\caption{\em Parameters for the equation \eqref{1-D_Graef}-\eqref{bc1}.}
\label{Graef_parameters}
\end{table}

\vspace{0.25cm}
\noindent\large{\textbf{Test Case 3.}}
We consider the fourth-order nonlinear equation describing the steady laminar flow of a viscous incompressible fluid in a porous channel \cite{Xu2010}.
For simplicity, we consider the one-dimensional problem that seeks $u$ such that
\begin{eqnarray}
\label{1D_channel_flows}
&&\frac{\td^4u}{\td x^4}+\gamma(x\frac{\td^3u}{\td x^3}+3\frac{\td^2u}{\td x^2})+R(u\frac{\td^3u}{\td x^3}-\frac{\td u}{\td x}\frac{\td^2u}{\td x^2})=0,\quad 0<x<1,
\\
&&
u(0)=0,\quad u''(0)=0,\quad u(1)=1,\quad u'(1)=0,\label{bcc}
\end{eqnarray}
where $R$ is the cross-flow Reynolds number and $\gamma$ is a physical constant related to the wall expansion ratio. Xu et al. \cite{Xu2010} have proven that the problem \eqref{1D_channel_flows}-\eqref{bcc} admits multiple solutions by analytic approaches. Three solutions were found by homotopy analysis method (HAM) in \cite{Liao2012} for the setting $R=-11$ and $\gamma=1.5$.

\begin{figure}
\centering
\includegraphics[scale=0.7]{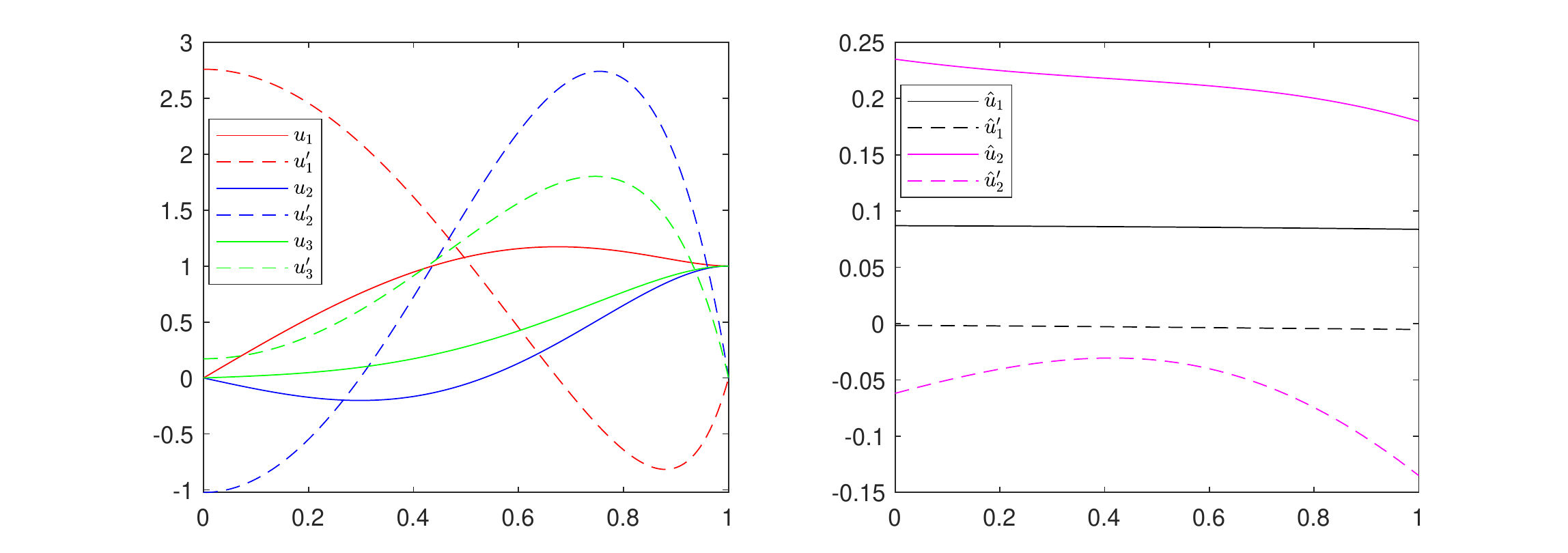}
\caption{\em Identified solutions and their derivatives of the channel flows equation \eqref{1D_channel_flows}-\eqref{bcc} by least square or neural network deflation. All correct solutions, $u_1$, $u_2$ and $u_3$, are identified with special networks for boundary conditions. Spurious solutions, $\bar{u}_1$ and $\bar{u}_2$, are found if the special networks are not used.}
\label{channel_flows_solutions}
\end{figure}

In our experiments, we take the same $R$ and $\gamma$ as in  \cite{Liao2012}. The special network for the boundary condition \eqref{bcc} is chosen as
\begin{equation}
u(x;\bm\theta,c)=x(x-1)^2(x^2\hat{u}(x;\bm\theta)+c)e^{2x}+\sin(\pi x/2),
\end{equation}
where $c$ is a network parameter to be trained together with $\bm{\theta}$. In this case, we initialize the bias of the third layer ${\bm b}_3={\bm 0}$ and $c\sim U[-5, 0]$. Other network parameters are initialized as mentioned above. Firstly, one solution $u_1$ is found by the least-square method \eqref{07}. Next, the second solution $u_2$ is obtained by the deflation \eqref{08} with deflation source $u_1$ ($p_1=2$). Moreover, the third solution $u_3$ is obtained by the deflation \eqref{08} with deflation sources $u_1$ and $u_2$ ($p_1=p_2=2$). Corresponding parameters are shown in Table \ref{channel_flows_parameters}. The three found solutions and their first derivatives are plotted in Figure \ref{channel_flows_solutions}, which are the same solutions found in \cite{Liao2012}.

Also, a comparison test is performed to seek $u_2$ by the deflation \eqref{16} with the same setting as above, except for using a generic solution network without special structures for boundary conditions. We find two solutions, denoted by $\hat{u}_1$ and $\hat{u}_2$, using $\lambda=1$ and $\lambda=100$. Neither of them takes the prescribed boundary value 0 at $x=0$ or 1 at $x=1$ and, hence, they are spurious solutions (see Figure \ref{channel_flows_solutions}).

\begin{table}
\centering
\begin{tabular}{|c|c|c|c|c|c|}
  \hline
  ~ & $u_1$ & $u_2$ & $u_3$ & $\hat{u}_1$ & $\hat{u}_2$\\\hline
  $N_\text{I}$ & 20000 & 10000 & 20000 & 10000 & 10000\\\hline
  $N_\text{p}$ & 1000 & 1000 & 1000 & 1000 & 1000\\\hline
  $I_\text{lr}$ & $[10^{-3},10^{-2}]$ & $[10^{-3},10^{-2}]$ & $[10^{-3},10^{-2}]$ & $[10^{-3},10^{-2}]$ & $[10^{-3},10^{-2}]$\\\hline
  deflation source & / & $u_1$ ($p_1=2$)& $u_1$, $u_2$  ($p_1=p_2=2$)& $u_1$ ($p_1=2$)& $u_1$ ($p_1=2$)\\\hline
  $\alpha$ & / & $1$ & $1$ & $1$ & $1$\\
  \hline
\end{tabular}
\caption{\em Parameters for the channel flows equation \eqref{1D_channel_flows}-\eqref{bcc}.}
\label{channel_flows_parameters}
\end{table}

\vspace{0.25cm}
\noindent\large{\textbf{Test Case 4.}}
We consider the following second-order problem that seeks $u$ such that
\begin{eqnarray}\label{1D_bootstrapping1}
&&\frac{\td^2u}{\td x^2}=f(u), \quad 0<x<1,\\
&& u'(0)=0,\quad u(1)=0,\label{bou1}
\end{eqnarray}
where $f(u)$ is a polynomial function of $u$. The existence of multiple solutions for the problem \eqref{1D_bootstrapping1} has been studied by the bootstrapping method \cite{Hao2014}.

\begin{figure}
\centering
\includegraphics[scale=0.6]{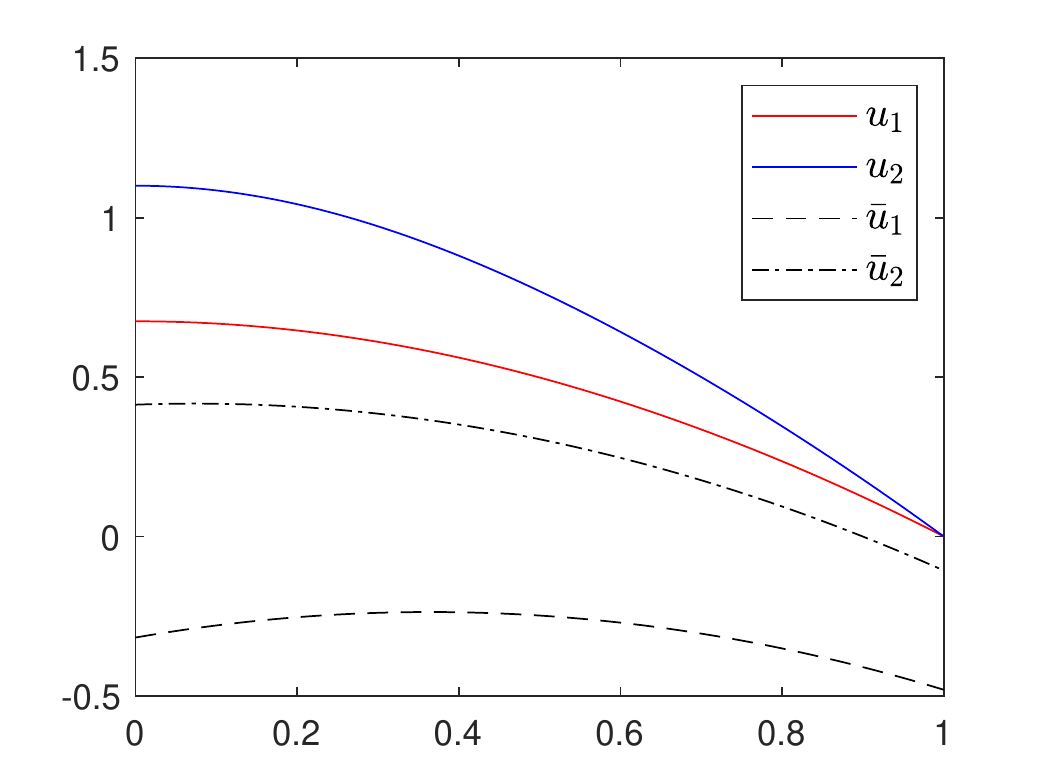}
\caption{\em Identified solutions of the nonlinear problem \eqref{1D_bootstrapping1}-\eqref{bou1} with $f(u)=\lambda(1+u^4)$ by least square or neural network deflation. All correct solutions, $u_1$ and $u_2$, are identified with special networks for boundary conditions. Spurious solutions, $\bar{u}_1$ and $\bar{u}_2$, are found if the special networks are not used.}
\label{1D_bootstrapping1_solutions}
\end{figure}

First, we set the right-hand side of the problem \eqref{1D_bootstrapping1} as $f(u)=\lambda(1+u^4)$. It is shown in \cite{Hao2014} that there are two solutions for $0<\lambda<\lambda^*=1.30107$. In our experiments, we take $\lambda=1.2$. The special network for the boundary condition \eqref{bou1} is given by
\begin{equation}\label{17}
u(x;\bm\theta)=x^2\hat{u}(x;\bm\theta)-\hat{u}(1;\bm\theta).
\end{equation}
The first solution $u_1$ is found by the least-square method \eqref{07} and the second solution $u_2$ is found by the deflation \eqref{08} with deflation source $u_1$ ($p_1=2$). Similarly to preceding cases, we perform a comparison test without the special network structure for boundary conditions and two spurious solutions $\hat{u}_1$ (for $\lambda=1$) and $\hat{u}_2$ (for $\lambda=100$) are found by the deflation \eqref{16}. The parameters for all these solutions are shown in Table \ref{1D_bootstrapping1_parameters} and all solutions are plotted in Figure \ref{1D_bootstrapping1_solutions}.

\begin{table}
\centering
\begin{tabular}{|c|c|c|c|c|}
  \hline
  ~ & $u_1$ & $u_2$ & $\hat{u}_1$ & $\hat{u}_2$ \\\hline
  $N_\text{I}$ & 10000 & 10000 & 10000 & 10000\\\hline
  $N_\text{p}$ & 1000 & 1000 & 1000 & 1000\\\hline
  $I_\text{lr}$ & $[10^{-3},10^{-2}]$ & $[10^{-3},10^{-2}]$ & $[10^{-3},10^{-2}]$ & $[10^{-3},10^{-2}]$\\\hline
  deflation source & / & $u_1$ ($p_1=2$) & $u_1$ ($p_1=2$)& $u_1$ ($p_1=2$)\\\hline
  $\alpha$ & / & $1$ & $1$ & $1$ \\
  \hline
\end{tabular}
\caption{\em Parameters for the nonlinear problem \eqref{1D_bootstrapping1}-\eqref{bou1} with $f(u)=\lambda(1+u^4)$.}
\label{1D_bootstrapping1_parameters}
\end{table}

\begin{figure}
\centering
\includegraphics[scale=0.6]{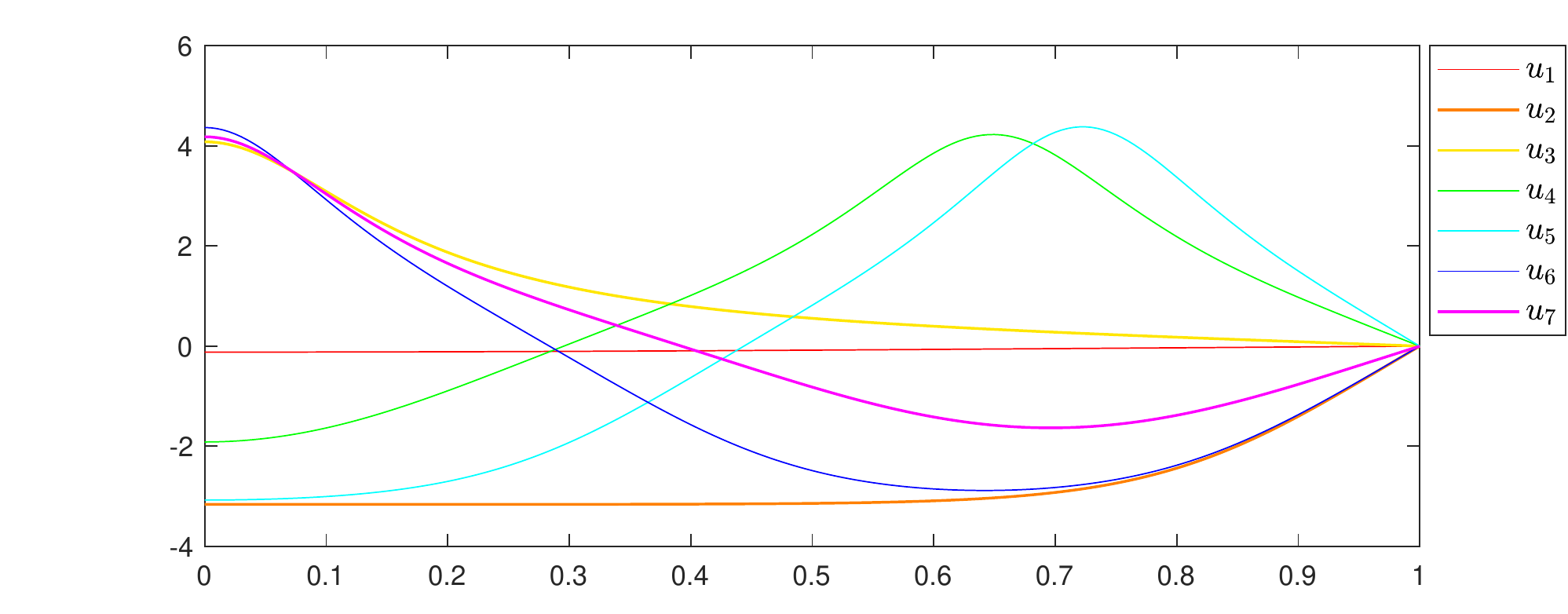}
\caption{\em Identified solutions of the nonlinear equation \eqref{1D_bootstrapping1}-\eqref{bou1} with $f(u)=-\frac{\pi^2}{4}u^2(u^2-10)$ by the deflation.}
\label{1D_bootstrapping2_solutions}
\end{figure}

Second, we repeat the test by choosing $f(u)=-\frac{\pi^2}{4}u^2(u^2-10)$. \cite{Hao2014} has proved that there exist eight solutions in total. Note that $u_0=0$ is a trivial solution. In this case, we start from the deflation \eqref{08} with the special network \eqref{17} and the deflation source $u_0$ ($p_0=2$) to find the first solution $u_1$, which is quite close to $u_0$. We would like to emphasize that it is sufficient to use the deflation without structure probing initializations to identify $u_0$ and $u_1$. However, we were not able to identify any other solutions without the structure probing initialization even if we tried our best to tune parameters and use different random initialization. To perform a wider search for other solutions, we employ the following structure probing initialization
\begin{equation}\label{19}
u_J(x;\bm\theta,c_j)=x^2\hat{u}(x;\bm\theta)-\hat{u}(1;\bm\theta)+\underset{j=1}{\overset{J}{\sum}}c_j\cos((2j-1)\pi x/2),
\end{equation}
with initial setting $c_j=0$ for $j=1,\cdots,J-1$ and $c_J\sim U(-5,5)$. Two solutions, denoted by $u_2$ and $u_3$, are found by the deflation \eqref{08} with source $u_0$ ($p_0=2$) and the structure probing network \eqref{19} with $J=1$. Another two solutions, denoted by $u_4$ and $u_6$, are found by the deflation \eqref{08} with source $u_0$ ($p_0=2$) and the network \eqref{19} with $J=2$. Two more solutions, denoted by $u_5$ and $u_7$, are found by the deflation \eqref{08} with deflation sources $u_4$ ($p_4=2$) and $u_6$ ($p_6=2$), respectively, and the network \eqref{19} with $J=2$. Corresponding parameters, including the initial value of $c_J$ actually randomized for each solution, are listed in Table \ref{1D_bootstrapping2_parameters}. All the $7$ nontrivial solutions are plotted in Figure \ref{1D_bootstrapping2_solutions}.

\begin{table}
\centering
\begin{tabular}{|c|c|c|c|c|}
  \hline
  ~ & $u_1$ & $u_2$ & $u_3$ & $u_4$ \\\hline
  $N_\text{I}$ & 5000 & 5000 & 10000 & 20000\\\hline
  $N_\text{p}$ & 1000 & 1000 & 1000 & 1000\\\hline
  $I_\text{lr}$ & $[10^{-3},10^{-2}]$ & $[10^{-3},10^{-2}]$ & $[10^{-4},10^{-3}]$ & $[10^{-4},10^{-3}]$\\\hline
  $J$ & / & 1 & 1 & 2\\\hline
  initial $c_J$ & / & $-3.48$ & $4.61$ & $-3.67$\\\hline
  deflation source & $u_0$ ($p_0=2$)& $u_0$ ($p_0=2$)& $u_0$ ($p_0=2$)& $u_0$ ($p_0=2$)\\\hline
  \hline
  ~ & $u_5$ & $u_6$ & $u_7$ & \\\hline
  $N_\text{I}$ & 20000 & 20000 & 20000 &\\\hline
  $N_\text{p}$ & 1000 & 1000 & 1000 &\\\hline
  $I_\text{lr}$ & $[10^{-4},10^{-3}]$ & $[10^{-4},10^{-3}]$ & $[10^{-4},10^{-3}]$ &\\\hline
  $J$ & 2 & 2 & 2 &\\\hline
  initial $c_J$ & $-4.12$ & $3.64$ & $3.44$ &\\\hline
  deflation source & $u_4$ ($p_4=2$) & $u_0$ ($p_0=2$) & $u_6$ ($p_6=2$) &\\\hline
\end{tabular}
\caption{\em Parameters for the nonlinear problem \eqref{1D_bootstrapping1}-\eqref{bou1} with $f(u)=-\frac{\pi^2}{4}u^2(u^2-10)$.}
\label{1D_bootstrapping2_parameters}
\end{table}

\subsection{Numerical tests in high-dimension} In this subsection, we will provide numerical tests in high-dimensional domains ($d \geq 2$).

\vspace{0.25cm}
\noindent\textbf{\large{Test Case 5.}}  We consider 2-D Yamabe's equation that seeks $u$ such that
\begin{equation}\label{2D_Yamabe}
\begin{split}
-8\Delta u-0.1u+\frac{u^5}{|\bm x|^3}=&0,\quad\text{in~}\Omega=\{\bm x\in \mathbb R^2: r<|\bm x|<R\},
\\
u=&1, \quad\text{on~}\partial\Omega,
\end{split}
\end{equation}
where $r$ and $R$ are set as 1 and 100. Nine solutions were found by using non-network deflation techniques and various initial guesses in \cite{Farrell2015}.

In our experiments, the solutions are approximated by the following special network
\begin{equation}\label{10}
u_J(\bm x;\bm \theta)=\hat{u}(\bm x;\bm \theta)\sin\left(\pi\frac{|\bm x|-r}{R-r}\right)+1
\end{equation}
if the random initialization without the structure probing technique is used,
or the following network
\begin{equation}\label{20}
u_J(\bm x;\bm \theta,c_j)=\hat{u}(\bm x;\bm \theta)\sin\left(\pi\frac{|\bm x|-r}{R-r}\right)+\underset{j=1}{\overset{J}{\sum}}c_j\sin(j\pi\frac{|\bm x|-r}{R-r})+1
\end{equation}
with the structure probing initialization, where the initial values are $c_j=0$ for $j=1,\cdots,J-1$ and $c_J\sim U(-1,1)$. Note that both \eqref{10} and \eqref{20} satisfy the given boundary condition automatically.

\begin{figure}
\centering
\includegraphics[scale=1.0]{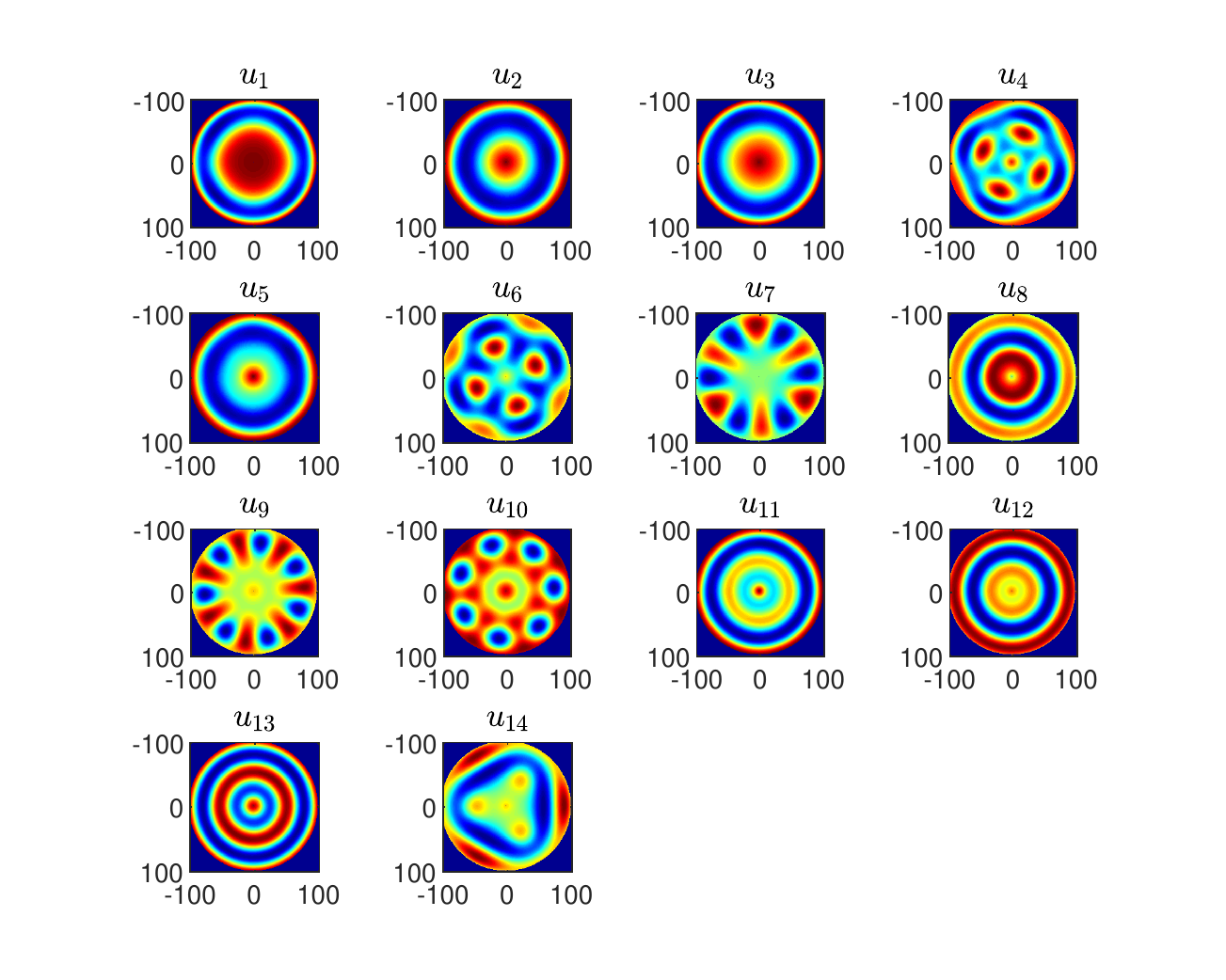}
\caption{\em Identified solutions of the 2-D Yamabe's equation \eqref{2D_Yamabe}.}
\label{2-D_Yamabe_solutions}
\end{figure}

In our proposed framework of the network-based structure probing deflation with a varying shift, we always follow the four steps: 1) use the least-square method \eqref{07} to find the first few solutions; 2) use neural network deflation without structure probing and varying shifts to find other solutions; 3) use structure probing deflation without varying shifts to find more distinct solutions; 4) finally, use structure probing deflation with varying shifts to find extra distinct solutions. Following these procedures, we find $14$ solutions in total for the 2-D Yamabe's equation as plotted in Figure \ref{2-D_Yamabe_solutions} with parameters specified in Table \ref{2-D_Yamabe_parameters}.

More precisely, $u_1$ and $u_{11}$ are found by the least-square method \eqref{07} and the others are found by the deflation \eqref{08} with previously found solutions as deflation sources ($p_k=2$ for all $k$). In deflation, we employ the technique of varying shifts in deflation operators, which helps to find more distinct solutions. All solutions are found by using networks \eqref{10} or \eqref{20} (specified in Table \ref{2-D_Yamabe_parameters}) with their corresponding initialization as mentioned previously, except that we take the network \eqref{10} with $2-u_9$ as the initial guess to find $u_{10}$. We would like to remark that both the structure probing initialization and the varying shifts are key techniques to find more distinct solutions for high-dimensional problems. Without any of them, we cannot find $14$ distinct solutions even if we have tried our best to tune parameters with commonly used random initialization in the literature.

\begin{table}
\centering
\begin{tabular}{|c|c|c|c|c|c|}
  \hline
  ~ & $u_1$ & $u_2$ & $u_3$ & $u_4$ & $u_5$\\\hline
  $N_\text{I}$ & 2000 & 2000 & 2000 & 5000 & 2000 \\\hline
  $N_\text{p}$ & 10000 & 10000 & 10000 & 10000 & 10000 \\\hline
  $I_\text{lr}$ & $[10^{-3},10^{-1}]$ & $[10^{-2},10^{-1}]$ & $[10^{-2},10^{-1}]$ & $[10^{-2},10^{-1}]$ & $[10^{-2},10^{-1}]$\\\hline
  network & \eqref{10} & \eqref{10} & \eqref{10} & \eqref{10} & \eqref{10}\\\hline
  $\alpha$ & / & 1 & 1 & $[0.01,100]$ & $[0.01,100]$\\\hline
  deflation source & / & $u_1$ & $u_2$ & $u_3$ & $u_1$\\\hline
  \hline
  ~ & $u_6$ & $u_7$ & $u_8$ & $u_9$ & $u_{10}$\\\hline
  $N_\text{I}$ & 5000 & 10000 & 20000 & 20000 & 10000 \\\hline
  $N_\text{p}$ & 10000 & 10000 & 10000 & 20000 & 10000 \\\hline
  $I_\text{lr}$ & $[10^{-2},10^{-1}]$ & $[10^{-2},10^{-1}]$ & $[10^{-2},10^{-1}]$ & $[10^{-2},10^{-1}]$ & $[10^{-2},10^{-1}]$\\\hline
  network & \eqref{10} & \eqref{10} & \eqref{10} & \eqref{10} & \eqref{10}\\\hline
  $\alpha$ & $[0.01,10]$ & $[0.01,10]$ & $[0.01,10]$ & $[0.01,10]$ & $[0.01,10]$\\\hline
  deflation source & $u_1$,$u_4$ & $u_1$,$u_2$ & $u_1$,$u_2$ & $u_1$,$u_2$ &$u_9$ \\\hline
  \hline
  ~ & $u_{11}$ & $u_{12}$ & $u_{13}$ & $u_{14}$ & \\\hline
  $N_\text{I}$ & 2000 & 10000 & 10000 & 10000 & \\\hline
  $N_\text{p}$ & 10000 & 10000 & 10000 & 10000 &  \\\hline
  $I_\text{lr}$ & $[10^{-3},10^{-1}]$ & $[10^{-2},10^{-1}]$ & $10^{-2}$ & $[10^{-2},10^{-1}]$ & \\\hline
  network & \eqref{20} ($J=4$) & \eqref{20} ($J=4$)& \eqref{20} ($J=4$)& \eqref{10} &\\\hline
  $\alpha$ & / & 0.01 & $[0.01,10]$ & 1 & \\\hline
  deflation source & / & $u_{11}$ & $u_8$,$u_{11}$ & $u_{11}$ &\\\hline
\end{tabular}
\caption{\em Parameters for the 2-D Yamabe's equation \eqref{2D_Yamabe} ($p_k=2$ for all deflation sources for the solutions obtained by the deflation).}
\label{2-D_Yamabe_parameters}
\end{table}

\vspace{0.25cm}
\noindent\textbf{\large{Test Case 6.}}   The high-dimensional Yamabe's equation seeks $u$ such that
\begin{equation}\label{HD_Yamabe}
\begin{split}
-\frac{4(d-1)}{(d-2)}\Delta u-0.125u+\frac{u^\frac{d+2}{(d-2)}}{|\bm x|^3}=&0,\quad\text{in~}\Omega=\{1<|\bm x|<100\},
\\
u=&1, \quad\text{on~}\partial\Omega,
\end{split}
\end{equation}
where $d\geq 3$ is the dimension of the problem.

We continue applying the network \eqref{10} without structure probing initialization and the network \eqref{20} with the structure probing initialization as solution networks to solve Yamabe's equation when $d=3$ and $d=6$. The initialization parameters are the same as in the 2-D case.

\begin{figure}
\centering
\includegraphics[scale=1.0]{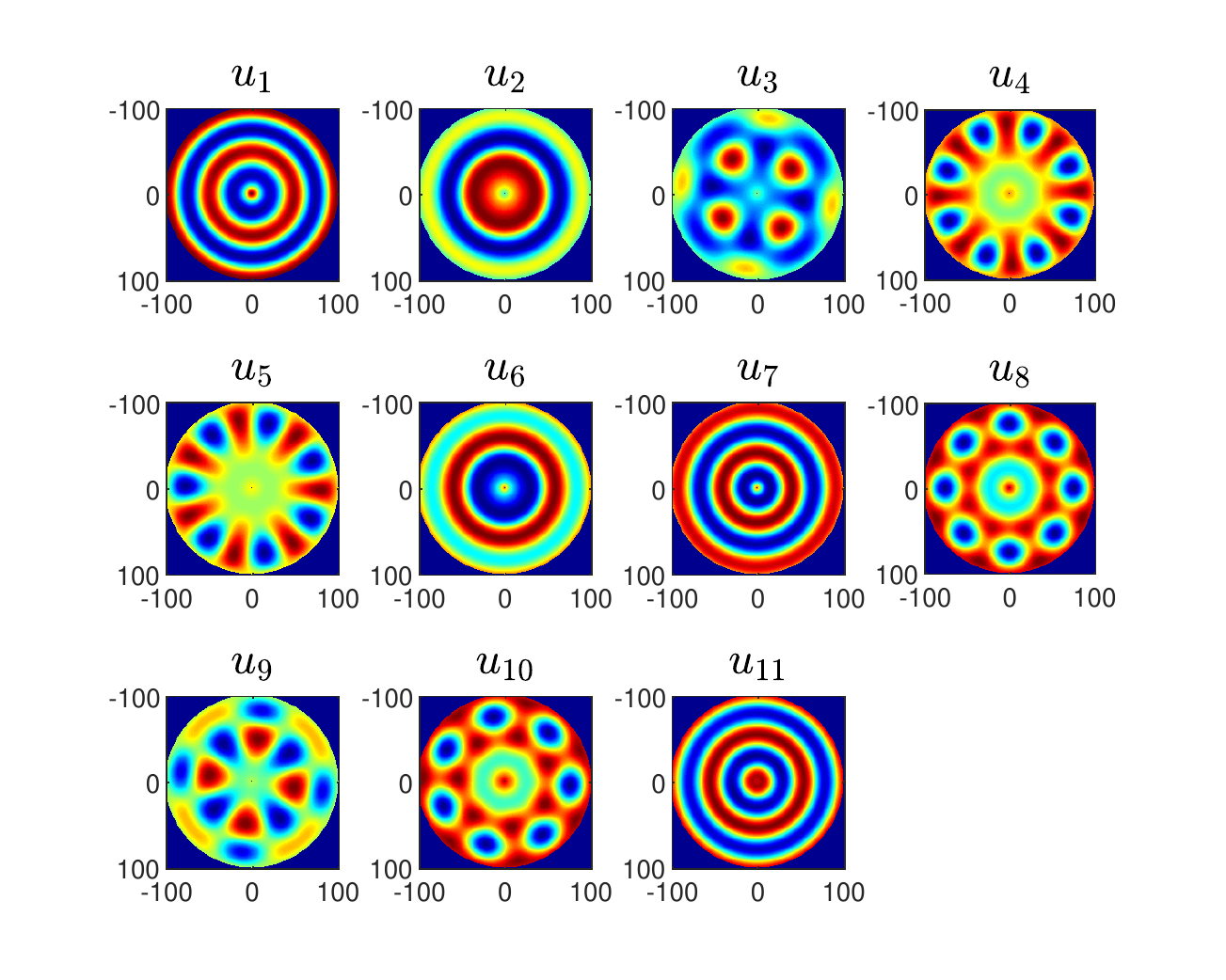}
\caption{\em Identified solutions of the 3-D Yamabe's equation \eqref{HD_Yamabe}. We visualize these solutions by projecting them in the first two coordinates.}
\label{3-D_Yamabe_solutions}
\end{figure}

\begin{figure}
\centering
\includegraphics[scale=1.0]{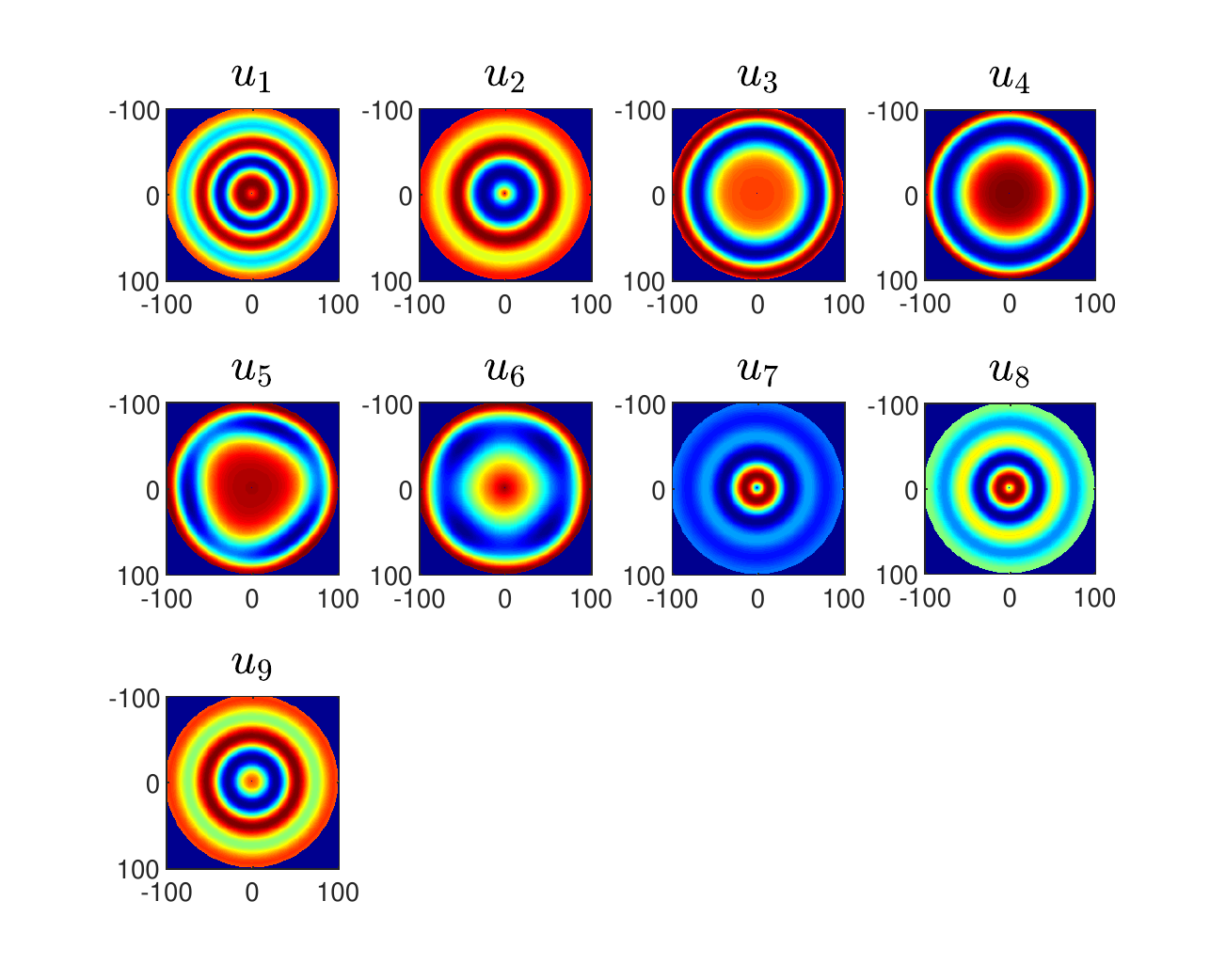}
\caption{\em Identified solutions of the 6-D Yamabe's equation \eqref{HD_Yamabe}. We visualize these solutions by projecting them in the first two coordinates.}
\label{6-D_Yamabe_solutions}
\end{figure}

Again, in our proposed framework of the network-based structure probing deflation with a varying shift, we follow the four steps: 1) use the least-square method \eqref{07} to find the first few solutions; 2) use the deflation without structure probing and varying shifts to find other solutions; 3) use structure probing deflation without varying shifts to find more distinct solutions; 4) finally, use structure probing deflation with varying shifts to find extra distinct solutions. Following these procedures, we obtain $11$ solutions when $d=3$ and $9$ solutions when $d=6$. The corresponding parameters are shown in Tables \ref{3-D_Yamabe_parameters} and \ref{6-D_Yamabe_parameters} for $d=3$ and $d=6$, respectively. The solutions are visualized in Figures \ref{3-D_Yamabe_solutions} and \ref{6-D_Yamabe_solutions} for $d=3$ and $d=6$, respectively. We would like to remark that both the structure probing initialization and the varying shifts are key techniques to find more distinct solutions for high-dimensional problems. Without any of them, we cannot find several distinct solutions even if we have tried our best to tune parameters with commonly used random initialization in the literature.

In these tests, the deflation powers $p_k$ are set as $2$ for all $k$ whenever deflation is used. In the case of $d=3$, most networks are initialized using \eqref{10} or \eqref{20}, except for $u_8$, $u_9$ and $u_{10}$, which are found by using initial guesses $2-u_4$, $2-u_3$ and $2-u_5$, respectively. In the case of $d=6$, we also try the initialization with a constant minus a known solution. However, this initialization method does not lead to new solutions.

\begin{table}
\centering
\begin{tabular}{|c|c|c|c|c|c|}
  \hline
  ~ & $u_1$ & $u_2$ & $u_3$ & $u_4$ & $u_5$\\\hline
  $N_\text{I}$ & 20000 & 20000 & 20000 & 20000 & 20000 \\\hline
  $N_\text{p}$ & 10000 & 10000 & 10000 & 10000 & 10000 \\\hline
  $I_\text{lr}$ & $[10^{-2},10^{-1}]$ & $[10^{-2},10^{-1}]$ & $[10^{-2},10^{-1}]$ & $[10^{-2},10^{-1}]$ & $[10^{-2},10^{-1}]$\\\hline
  network & \eqref{10} & \eqref{10} & \eqref{10} & \eqref{10} & \eqref{10}\\\hline
  $\alpha$ & / & $[0.01,10]$ & 1 & 0.1 & 0.01\\\hline
  deflation source & / & $u_1$ & $u_1$,$u_2$ & $u_1$,$u_2$ & $u_1$,$u_2$\\\hline
  \hline
  ~ & $u_6$ & $u_7$ & $u_8$ & $u_9$ & $u_{10}$\\\hline
  $N_\text{I}$ & 20000 & 20000 & 20000 & 20000 & 20000 \\\hline
  $N_\text{p}$ & 10000 & 10000 & 10000 & 10000 & 10000 \\\hline
  $I_\text{lr}$ & $[10^{-2},10^{-1}]$ & $[10^{-2},10^{-1}]$ & $[10^{-2},10^{-1}]$ & $[10^{-2},10^{-1}]$ & $[10^{-2},10^{-1}]$\\\hline
  network & \eqref{20} ($J=4$) & \eqref{20} ($J=6$) & \eqref{20} ($J=4$) & \eqref{10} & \eqref{10}\\\hline
  $\alpha$ & 0.01 & 0.1 & $[0.01,10]$ & $[0.01,10]$ & $[0.01,10]$\\\hline
  deflation source & $u_1$,$u_2$ & $u_1$,$u_2$ & $u_1$,$u_4$ & $u_1$,$u_2$,$u_3$ &$u_1$,$u_2$,$u_5$ \\\hline
  \hline
  ~ & $u_{11}$ &&&&\\\hline
  $N$ & 100  &  &  &  & \\\hline
  $N_\text{I}$ & 20000 &&&& \\\hline
  $N_\text{p}$ & 10000 &&&& \\\hline
  $I_\text{lr}$ & $[10^{-2},10^{-1}]$ &&&&\\\hline
  network & \eqref{20} ($J=4$) &&&&\\\hline
  $\alpha$ & $[0.01,10]$ &&&&\\\hline
  deflation source & $u_1$,$u_2$,$u_6$ &&&&\\\hline
\end{tabular}
\caption{\em Parameters for the 3-D Yamabe's equation \eqref{HD_Yamabe} ($p_k=2$ for all deflation sources for the solutions obtained by the deflation).}
\label{3-D_Yamabe_parameters}
\end{table}

\begin{table}
\centering
\begin{tabular}{|c|c|c|c|c|c|}
  \hline
  ~ & $u_1$ & $u_2$ & $u_3$ & $u_4$ & $u_5$\\\hline
  $N_\text{I}$ & 20000 & 20000 & 20000 & 20000 & 20000 \\\hline
  $N_\text{p}$ & 10000 & 10000 & 10000 & 10000 & 10000 \\\hline
  $I_\text{lr}$ & $[10^{-2},10^{-1}]$ & $[10^{-2},10^{-1}]$ & $[10^{-2},10^{-1}]$ & $[10^{-3},10^{-1}]$ & $[10^{-3},10^{-1}]$\\\hline
  network & \eqref{10} & \eqref{10} & \eqref{10} & \eqref{10} & \eqref{10}\\\hline
  $\alpha$ & / & $[0.01,10]$ & 0.1 & 10 & $[0.01,10]$\\\hline
  deflation source & / & $u_1$ & $u_1$ & $u_1$ & $u_1$,$u_2$\\\hline
  \hline
  ~ & $u_6$ & $u_7$ & $u_8$ & $u_9$ & \\\hline
  $N_\text{I}$ & 20000 & 20000 & 20000 & 20000 & \\\hline
  $N_\text{p}$ & 10000 & 10000 & 10000 & 10000 & \\\hline
  $I_\text{lr}$ & $[10^{-3},10^{-2}]$ & $[10^{-3},10^{-2}]$ & $[10^{-2},10^{-1}]$ & $[10^{-2},10^{-1}]$ & \\\hline
  network & \eqref{10} & \eqref{20} ($J=6$) & \eqref{20} ($J=6$) & \eqref{20} ($J=6$) & \\\hline
  $\alpha$ & $[0.1,1]$ & / & $[0.01,10]$ & 1 & \\\hline
  deflation source & $u_1$,$u_2$ & / & $u_7$ & $u_7$ & \\\hline
\end{tabular}
\caption{\em Parameters for the 6-D Yamabe's equation \eqref{HD_Yamabe} ($p_k=2$ for all deflation sources for the solutions obtained by the deflation).}
\label{6-D_Yamabe_parameters}
\end{table}

{
\vspace{0.25cm}
\noindent\textbf{\large{Test Case 7.}} In the last example, we consider the following reaction-diffusion system applied in the modeling of the chemical reaction with two components \cite{Nicolis1977} and irregular patterns \cite{Pearson1993},
\begin{equation}\label{3-D_chemical}
\begin{cases}
\mathcal{D}_1(u,v):=\varepsilon_u\Delta u-uv^2+F(1-u)=0\\
\mathcal{D}_2(u,v):=\varepsilon_v\Delta v+uv^2-(F+k)v=0
\end{cases}\quad\text{in}~\Omega,
\end{equation}
with Dirichlet boundary conditions
\begin{equation}
u=1\text{ and }v=0~\text{on}~\partial\Omega.
\end{equation}
In this case, $\Omega$ is set as a more complicated domain in $\mathbb{R}^3$ formulated by
\begin{equation}
\Omega=\{{\bm x}\in\mathbb{R}^3:|{\bm x}|<\rho({\bm x}):=1+0.1\sin(5\theta(x_1+\text{i}x_2))\},
\end{equation}
where $\theta(z)$ means the argument of a complex number $z$. See Fig. \ref{3-D_chemical_domain} for the visualization of $\Omega$. Note the system \eqref{3-D_chemical} has a pair of trivial solutions $u_0\equiv1$ and $v_0\equiv0$.

\begin{figure}
\centering
\includegraphics[scale=0.5]{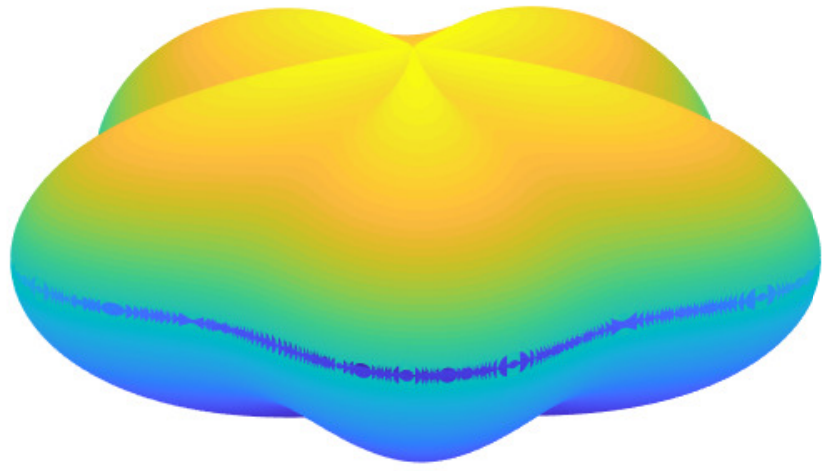}
\caption{\em The problem domain of the 3-D reaction-diffusion system.}
\label{3-D_chemical_domain}
\end{figure}

In \cite{Wang2018}, the authors solve the problem \eqref{3-D_chemical} in a 2-D square by a spectral collocation method, obtaining a vast number of solutions with residuals less than $10^{-9}$. However, it is quite challenging to solve the problem \eqref{3-D_chemical} in a 3-D complicated domain by most conventional approaches (e.g., FDM and spectral methods).

Our network-based strategy is to construct two special networks $u(\bm x;\bm{\theta}_u)$ and $v(\bm x;\bm{\theta}_v)$ to approximate $u$ and $v$, respectively. Specifically, we let
\begin{eqnarray}
&&u(\bm x;\bm{\theta}_u)=\hat{u}(\bm x;\bm{\theta}_u)(|{\bm x}|^2-\rho^2({\bm x}))+1,\\
&&v(\bm x;\bm{\theta}_v)=\hat{v}(\bm x;\bm{\theta}_v)(|{\bm x}|^2-\rho^2({\bm x})),
\end{eqnarray}
which automatically satisfy $u(\bm x;\bm{\theta}_u)=1$ and $v(\bm x;\bm{\theta}_v)=0$ on $\partial\Omega$. If we use the original least squares method in \eqref{03}, only the trivial solutions can be found. Therefore, we train the networks by the following deflation
\begin{multline}\label{23}
\underset{\bm{\theta}_u,\bm{\theta}_v}{\min}~L_\text{ND}(\bm{\theta}_u,\bm{\theta}_v):=\left(\overset{K}{\underset{k=1}{\sum}}\left(\|u(\bm x;\bm{\theta}_u)-u_k(\bm x)\|_{L^2(\Omega)}^{-p_k}+\|v(\bm x;\bm{\theta}_v)-v_k(\bm x)\|_{L^2(\Omega)}^{-p_k}\right)+\alpha\right)\\
\cdot\left(\|\mathcal{D}_1(u(\bm x;\bm{\theta}_u),v(\bm x;\bm{\theta}_v))\|_{L^2(\Omega)}^2+\|\mathcal{D}_2(u(\bm x;\bm{\theta}_u),v(\bm x;\bm{\theta}_v))\|_{L^2(\Omega)}^2\right),
\end{multline}
where $\left\{u_k(\bm x),v_k(\bm x)\right\}_{k=1}^K$ are $K$ pairs of solutions that have already been obtained. We start the search by taking the trivial solutions $u_0$ and $v_0$ as deflation sources, and then take identified solutions as new deflation sources for the next search. Hyper-parameters are set as $N_\text{I}=10000$, $N_\text{p}=10000$, and $I_\text{lr}=[10^{-5},10^{-2}]$. Besides, we use a varying $\alpha$ with a range $[10^{-2},1]$. Deflation powers are set as $p_k=2$ for all sources. Finally, we find more than 100 distinct solutions, some of which are shown in Fig. \ref{3-D_chemical_solutions}. The residual errors of all identified solutions for Equation \eqref{3-D_chemical} are below $1.0\times10^{-3}$ and the corresponding values of the loss function in \eqref{23} are below $0.5\times10^{-3}$. 

\begin{figure}
\centering
\includegraphics[scale=1.0]{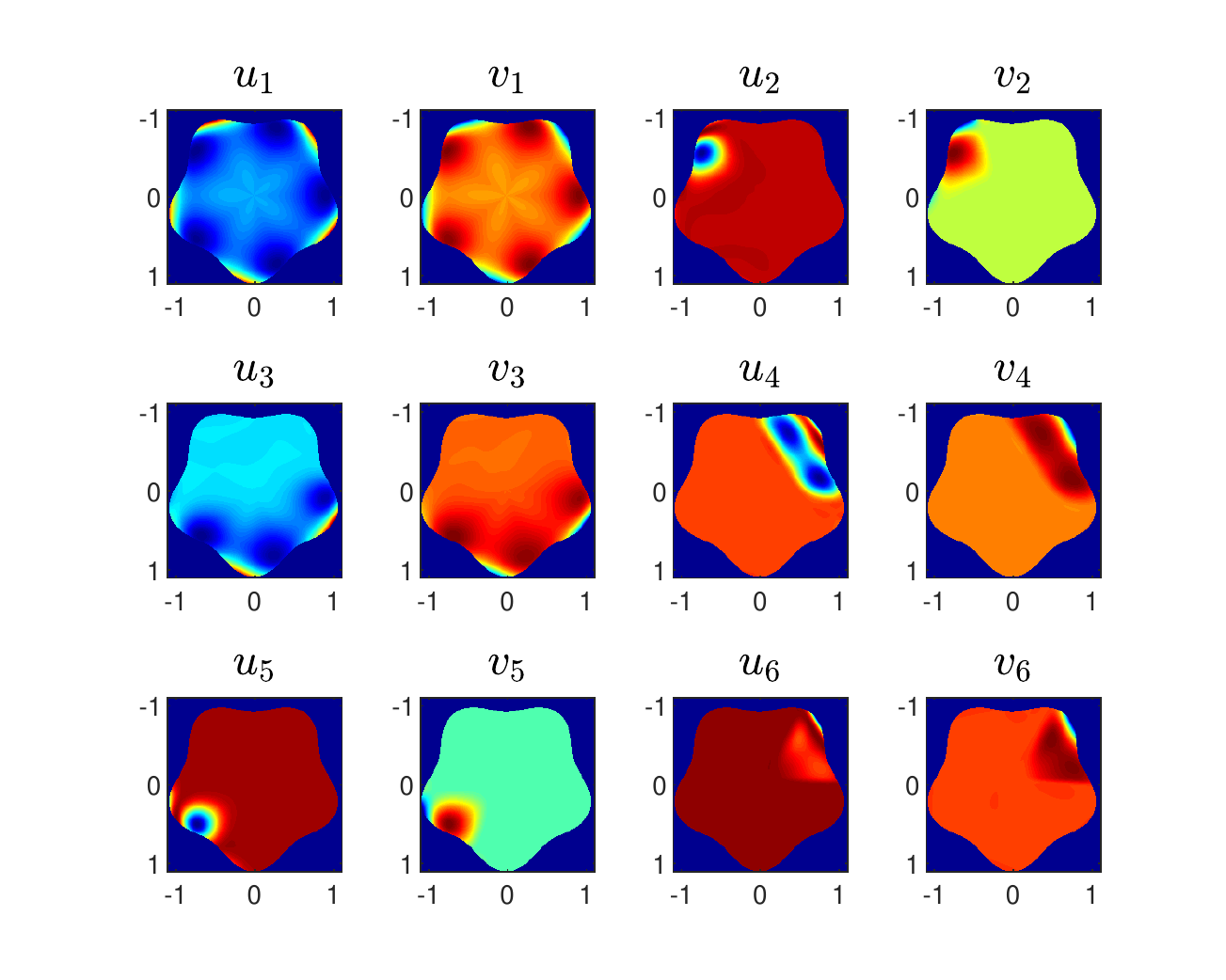}
\caption{\em Selected solution pairs $(u,v)$ of the 3-D reaction-diffusion system \eqref{3-D_chemical}. We visualize these solutions by projecting them in the first two coordinates.}
\label{3-D_chemical_solutions}
\end{figure}
}

\section{Conclusion}\label{sec:conclusion}
In this paper, we proposed the structure probing neural network deflation to find distinct solutions to nonlinear differential equations. The original optimization energy landscape of network-based methods is regularized by neural network deflation so that known solutions are no longer local minimizers while preserving unknown solutions as local minimizers. To obtain a new solution with the desired features, a structure probing algorithm is applied to obtain an initial guess that is close to the desired solution. Finally, special network structures that satisfy various boundary conditions automatically are introduced to simplify the objective function of network-based methods. These techniques form a new framework for identifying distinct solutions of nonlinear differential equations. Compared to existing methods, the proposed neural network deflation is capable of solving high-dimensional problems on complex domains with a lower computational cost and can identify more distinct solutions. { As a neural network-based PDE solver, structure probing neural network deflation may not provide highly accurate solutions. But these solutions are usually accurate enough for industrial applications and serve as a good initial guess for conventional methods as in \cite{Huang2019} to obtain highly accurate solutions efficiently.

Structure probing neural network deflation relies on the deflation operator proposed in \cite{Farrell2015} based on conventional discretization methods. Although the application of neural networks has conquered some disadvantages of the conventional deflation method, e.g., we can solve high-dimensional problems on complex domains and identify more solutions, the proposed method in this paper still inherits some disadvantages of the conventional deflation method. For example, when two solutions are very close to each other, the optimization landscape of the deflated loss using one solution as the deflation source becomes very steep at the other solution, making it very challenging to identify another solution. As in the conventional deflation method, it is crucial to choose appropriate powers $p_k$ for deflation sources as shown in our numerical tests. However, the parameter selection is still heuristic and problem-dependent. Learning how to choose parameters automatically is an important future direction. Network-based methods in general might need extra effort to deal with boundary conditions, which is not an issue of conventional methods. Designing more advanced optimization algorithms for constrained optimization in network-based methods would also be interesting in the future.
}

{\bf Acknowledgments.} Y. G. was partially supported by the Ministry of Education in Singapore under the grant MOE2018-T2-2-147. C. W. is partially supported by the US National Science Foundation Award DMS-1849483. H. Y. was partially supported by the US National Science Foundation under award DMS-1945029.

\bibliography{expbib}
\bibliographystyle{plain}

\end{document}